\documentclass[12pt]{article}
\usepackage{amsthm}
\usepackage{amsmath}
\usepackage{ bbold }
\usepackage[colorlinks,citecolor=blue,urlcolor=blue,filecolor=blue,backref=page]{hyperref}
\usepackage{graphicx}
\usepackage[english]{babel}
\usepackage[utf8]{inputenc}
\usepackage{amsmath}
\usepackage{amssymb}
\usepackage{graphicx}
\usepackage{hyperref}
\usepackage{graphicx}
\usepackage{comment}
\usepackage{mathrsfs}
\usepackage{amsthm}
\usepackage{bbm}
\usepackage{rotating}
\usepackage{lscape}
\usepackage{subfigure}
\usepackage[colorinlistoftodos]{todonotes}
\usepackage{lipsum}
\usepackage{comment}
\usepackage[round]{natbib}
\usepackage[colorlinks,citecolor=blue,urlcolor=blue,filecolor=blue,backref=page]{hyperref}
\usepackage{graphicx}
\newtheorem{theorem}{Theorem}

\newtheorem{plan}[theorem]{Plan}
\newtheorem{claim}[theorem]{Claim}
\newtheorem{remark}[theorem]{Remark}
\numberwithin{equation}{section}
\theoremstyle{plain}

\newcommand{\ignore}[1]{}

\newcommand*{\tran}{^{\mkern-1.5mu\mathsf{T}}}

\newcommand{\sF}{\mathcal{F}}

\newcommand{\IE}{\mathbb{E}}

\newcommand{\IR}{\mathbb{R}}

\newcommand{\IP}{\mathbb{P}}

\def\argmin{\mathop{\rm argmin}}

\newcommand{\blind}{0}

\begin{document}

\def\spacingset#1{\renewcommand{\baselinestretch}%
{#1}\small\normalsize} \spacingset{1}

\if0\blind
{
  \title{\bf {Bayesian modeling of time-varying conditional heteroscedasticity}}
  \author{Sayar Karmakar, 
       Arkaprava Roy\\
       University of Florida}
  \maketitle
} \fi


\if1\blind
{
  \bigskip
  \bigskip
  \bigskip
  \begin{center}
    {\LARGE\bf }
\end{center}
  \medskip
} \fi

\bigskip

\begin{abstract}
Conditional heteroscedastic (CH) models are routinely used to analyze financial datasets. The classical models such as ARCH-GARCH with time-invariant coefficients are often inadequate to describe frequent changes over time due to market variability. However, we can achieve significantly better insight by considering the time-varying analogs of these models. In this paper, we propose a Bayesian approach to the estimation of such models and develop a computationally efficient MCMC algorithm based on Hamiltonian Monte Carlo (HMC) sampling. We also established posterior contraction rates with increasing sample size in terms of the average Hellinger metric. The performance of our method is compared with frequentist estimates and estimates from the time constant analogs. To conclude the paper we obtain time-varying parameter estimates for some popular Forex (currency conversion rate) and stock market datasets. 
\end{abstract}

\textbf{Keywords:} Autoregressive model, B-splines, Hamiltonian Monte Carlo (HMC), Non-stationary, Posterior contraction,  Volatility

\section{Introduction}

For datasets observed over a long period, stationarity turns out to be an oversimplified assumption that ignores systematic deviations of parameters from constancy. Thus time-varying parameter models have been studied extensively in the literature of statistics, economics, and related fields. For example, the financial datasets, due to external factors such as war, terrorist attacks, economic crisis, political events, etc. exhibit deviation from time-constant stationary models. Accounting for such changes is crucial as otherwise time-constant models can lead to incorrect policy implications as pointed out by \citet{bai97}. Thus functional estimation of unknown parameter curves using time-varying models has become an important research topic today. In this paper, we analyze popular conditional heteroscedastic models such as AutoRegressive Conditional Heteroscedasticity (ARCH) and Generalized ARCH (GARCH) from a Bayesian perspective in a time-varying setup.
Before discussing our new contributions in this paper, we provide a brief overview of some previous works in this areas.

In the regression regime, time-varying models have garnered a lot of recent attention; see, for instance, \citet{fan99}, \citet{fan2000}, \citet{hoover98}, \citet{huang04}, \citet{lin01}, \citet{ramsay05}, \citet{zhang02} among others.  The models show time-heterogeneous relationship between response and predictors. Consider the following two regression models
$$\text{Model I: } y_i = x_i\tran \theta_i +e_i, \quad \text{Model II: } y_i = x_i\tran\theta_0 + e_i, \quad\quad i = 1, \ldots ,n,$$

\noindent where $x_i \in \IR^d$ ($i = 1,\ldots,n$) are the covariates, $\tran$ is the transpose, $\theta_0$ and $\theta_i = \theta(i/n)$ are the regression coefficients. Here, $\theta_0$ is a constant parameter and $\theta : [0, 1] \to  \mathbb{R}^d$ is a smooth function. Estimation of $\theta(\cdot)$ has been considered by \citet{hoover98}, \citet{cai07}) and \citet{zhouwu10} among others. One popular way to decide if there is an evidence to favor time-varying models over the time-constant analogue is to perform hypothesis testing. See, for instance, \citet{regression2012}, \citet{regression2015}, \citet{chow60}, \citet{brown75}, \citet{nabeya88},  \citet{leybourne89}, \citet{nyblom89}, \citet{ploberger89}, \citet{andrews93} and \citet{lin99}. \citet{zhouwu10} discussed obtaining  simultaneous confidence bands (SCB) in model I, i.e. with additive errors. However their treatment is heavily based on the closed-form solution and it does not extend to processes defined by a more general recursion.


For time-varying AR, MA, or ARMA processes, the results from time-varying linear regression can be naturally extended. However, such an extension is not obvious for conditional heteroscedastic (CH hereafter ) models. These are, by the simple definition of evolution is difficult to estimate even in the time-constant case. However, one cannot possibly ignore its usefulness in analyzing and predicting financial datasets. These models (even the simple time-constant ones) have remained primary tools for analyzing and forecasting certain trends for stock market datasets since \citet{engle82} introduced the classical ARCH model and \citet{bollerslev} extended it to a more general GARCH model. However, with the rapid dynamics of market vulnerability, the simple classical time-constant models fail in terms of estimation or prediction due to over-compensating the past data. Several references point out the necessity of extending these classical models to a set-up where the parameters can change across time, for example \citet{stuaricua2005}, \citet{engle2005spline} and \citet{fry08}. Consider the simple tvARCH(1) model
$$X_i= \sigma_i\zeta_i, \zeta_i \sim N(0,1), \sigma_i^2= \mu_0(i/n)+a_1(i/n)X_{i-1}^2.$$ Similar models can be defined for tvGARCH(1,1) as well where $\sigma_i^2$ has an additional recursive term involving $\sigma_{i-1}^2$
$$X_i= \sigma_i\zeta_i, \zeta_i \sim N(0,1), \sigma_i^2= \mu_0(i/n)+a_1(i/n)X_{i-1}^2+b_1(i/n)\sigma_{i-1}^2.$$
When the two recursive parameters in a GARCH model sum up to 1, i.e. $a_1+b_1=1$ it is usually called an integrated GARCH (iGARCH; or bubble garch/explosive garch by some authors) process which employing the above display can also be extended towards a time-varying analog i.e. $b_1(\cdot)=1-a_1(\cdot)$. A wide range of financial datasets exhibits iGARCH phenomena. 

In the parlance of time-varying parameter models in the CH setting, numerous works discussed the CUSUM-type procedure, for instance, \citet{kim00} for testing change in parameters of GARCH(1,1). \citet{kulperger05} studied the high moment partial sum process based on residuals and applied it to residual CUSUM tests in GARCH models. Interested readers can find some more change–point detection results in the context of CH models in \citet{chu95}, \citet{chen97}, \citet{lin1999testing}, \citet{kokoszka00} or \citet{andreou06}.

A time-varying framework and a pointwise curve estimation using M-estimators for locally stationary ARCH models were provided by \citet{dahlhaussubbarao2006}. Since then, while several pointwise approaches were discussed in the tvARMA and tvARCH case (cf. \citet{dahlhaus2009}, \citet{dahlhaussubbarao2006}, \citet{fry08}), pointwise theoretical results for estimation in tvGARCH processes were discussed in \citet{tvgarch2013} and \citet{rohan13} for GARCH(1,1) and GARCH($p$,$q$) models respectively. In a series of recent works \cite{sayar2020, sayarthesis} such models were discussed in wide generality. However, the focus remained frequentist, and the main goal accomplished there was to build simultaneous inference. One strong criticism for the CH type models remained that one needs a relatively large sample size ($n \sim 2000$) to achieve nominal coverage levels. The recursive definition of the models and a subsequent kernel-based method of estimating make it difficult to achieve satisfying results for relatively smaller sample sizes. This motivated us to explore a Bayesian way of building and estimating these models and use the posteriors to construct posterior estimates of the coefficient curves $\theta(\cdot)$.

In this paper, we develop a Bayesian estimation method for time-varying analogs of ARCH, GARCH, and iGARCH models. We model the time-varying functional parameters using cubic B-splines. In the context of general varying-coefficient modeling, spline bases are a popular choice for its convenience and flexibility \citep{hastie1993varying, gu1993smoothing, cai2000functional,biller2001Bayesian, huang2002varying, huang2004functional,amorim2008regression,fan2008statistical,yue2014Bayesian,franco2019unified}.  Specific to the literature of time-varying volatility modeling, B-spline-based models have also been explored \citep{engle2008spline,audrino2009splines,liu2016spline}.

Our contributions in this paper are two-fold. Towards the methodological development, note that the tvARCH, tvGARCH, and tviGARCH models require complex shape constraints on the coefficient functions. We achieve those by imposing different hierarchical structures on B-spline coefficients. The constraints are designed to be able to develop an efficient sampling algorithm based on gradient-based Hamiltonian Monte Carlo (HMC) \citep{neal2011mcmc, betancourt2015hamiltonian, betancourt2017conceptual,livingstone2019geometric}. Strong motivation towards implementing such a Bayesian methodology was to circumvent the requirement of a huge sample size which is almost essential for effective estimation using the frequentist and kernel-based methods. This requirement on sample size has been frequently pointed out in the literature of ARCH/GARCH models and thus this was one of our main motivations to see if a reasonable estimation scheme can be designed in a Bayesian way.

Secondly, the existing literature on obtaining posterior concentration rates for dependent data is thin, even for an extremely simple model. To the best of our knowledge, ours is the first such attempt towards a theoretical development for these models under Gaussian-link. Posterior contraction rates for these models with respect to the average Hellinger metric are established. The main challenge therein is to construct exponentially consistent tests for these classes of models. Using some recently developed tools from \cite{jeong2019frequentist, ning2020Bayesian} we have developed such tests. We first establish posterior contraction rates with respect to average log-affinity and then the same rate is transferred to the average Hellinger metric. The frequentist literature on inference about time-varying needs very stringent moment assumption and local stationarity assumptions which are often difficult to verify. Moreover, for econometric datasets, the existence of even the fourth moment is often questionable. Thus this paper offers some alternative way to estimate coefficients under lesser assumptions.

The rest of the paper is organized as follows. Section \ref{sec:model} describes the proposed Bayesian model in detail. Section~\ref{sec:comp} discusses an efficient computational scheme for the proposed method. We calculate posterior contraction rate in Section~\ref{sec:theo}. In Section~\ref{sec:sim} we study the performance of our proposed method in the light of. Section~\ref{sec:application} deals with real data application of the proposed methods for the three separate models and concludes with a brief interpretation of the results. We wrap the paper up with discussions, some concluding remarks, and possible future directions in Section~\ref{sec:discussion}. The supplementary materials contain theoretical proofs and some additional results.

\section{Modeling}
\label{sec:model}
We elaborate on the models and our Bayesian framework for time-varying analogs of three specific cases that are popularly used to analyze econometric datasets.
\subsection{tvARCH Model}\label{ssc:tvarchmodel}
Let $\{X_i\}$ satisfy the following time-varying ARCH($p$) model for $X_i$ given $\mathcal{F}_{i-1}=\{X_j: j\leq (i-1)\}$,
\begin{align}
    &X_i|\mathcal{F}_{i-1}\sim \mathrm{N}(0,\sigma_i^2) \\
    &\sigma_i^2=\mu(i/n) + \sum_{k=1}^p a_k(i/n) X_{i-k}^2\label{modelarch}
\end{align}
where the parameter functions $\mu(\cdot), a_i(\cdot)$ satisfy
\begin{align}\label{eq:parcondition}
    \mathcal{P}=\{\mu, a_k:\mu(x)\geq 0, 0\leq a_{k}(x)\leq 1, \sup_{x}\sum_{k}a_{k}(x)<1\}
\end{align}

\noindent In a Bayesian regime we put priors on $\mu(\cdot)$ and $a_i(\cdot)$. To respect the shape-constraints as imposed by $\mathcal{P}$ we reformulate the problem. With $B_j$ as the B-spline basis functions, let

\begin{align*}
\mu(x) =&\sum_{j=1}^{K_1}\exp(\beta_j)B_j(x)\\
     a_{k}(x)=&\sum_{j=1}^{K_2}\theta_{kj}M_{k}B_j(x), \quad 0\leq\theta_{kj}\leq 1,\\
    M_i=&\frac{\exp(\delta_i)}{\sum_{k=0}^p{\exp(\delta_k)}}, \quad i=1,\ldots,p,\\
    \delta_l\sim&N(0, c_1),\textrm{ for }0\leq l\leq p,\\
    \beta_{j}\sim & N(0, c_2)\textrm{ for } 1\leq j\leq K_1,\\
    \theta_{kj}\sim& U(0,1)\textrm{ for }1\leq k\leq p, 1\leq j\leq K_2.
\end{align*}

The prior induced by above construction is $\mathcal{P}$-supported. The verification is very straightforward. In above construction, $\sum_{j=0}^PM_j=1$. Thus $\sum_{j=1}^PM_j\leq 1$. Since $0\leq\theta_{kj}\leq 1$, $\sup_{x}a_i(x)\leq M_i$. Thus $\sup_x\sum_{i=1}^Pa_{i}(x)\leq \sum_{i=1}^PM_i\leq 1$. We have $\sum_{j=1}^PM_j\leq 1$ if and only if $\delta_0=-\infty$, which has probability zero. On the other hand, we also have $\mu(\cdot)\geq 0$ as we have $\exp(\beta_j)\geq 0$. Thus, the induced priors, described above are well supported in $\mathcal{P}$.

\subsection{tvGARCH Model}\label{ssc:tvgarchmodel}
Let $\{X_i\}$ satisfy the following time-varying GARCH($p$,$q)$ model for $X_i$ given $\mathcal{F}_{i-1}=\{X_j: j\leq (i-1)\}$,
\begin{align}
    &X_i|\mathcal{F}_{i-1}\sim \mathrm{N}(0,\sigma_i^2) \nonumber\\
    &\sigma_i^2=\mu(i/n) + \sum_{k=1}^p a_k(i/n) X_{i-k}^2+\sum_{j=1}^q b_j(i/n) \sigma^2_{i-j}\label{modelgarch}
\end{align}
Additionally we impose the following constraints on parameter space for the time-varying parameters, 
\begin{align}\label{eq:parconditiongarch}
    \mathcal{P}_1=\{\mu, a_i:\mu(x)\geq 0, 0\leq a_{i}(x), 0 \leq b_j(x), \sup_{x}\sum_{k}a_{k}(x)+\sum_{j}b_{j}(x)<1\}
\end{align}
The condition on the AR parameters imposed by (\ref{eq:parconditiongarch}) is somewhat popular in time-varying AR literature. See \cite{dahlhaussubbarao2006, subbarao2008, sayar2020} for example. Different from these references,  we additionally do not assume existence of any unobserved local-stationary process that are close to the observed process.

To proceed with Bayesian computation, we again put priors on the unknown functions $\mu(\cdot), a_i(\cdot)$ and $b_j(\cdot)$'s such that they are supported in $\mathcal{P}_1$. Again the restrictions imposed by (\ref{eq:parconditiongarch}) are respected. The complete description of prior is

\begin{align*}
\mu(x) =&\sum_{j=1}^{K_1}\exp(\beta_j)B_j(x)\\
     a_{k}(x)=&\sum_{j=1}^{K_2}\theta_{kj}M_{k}B_j(x), \quad 0\leq\theta_{kj}\leq 1,1\leq k\leq p,\\
     b_{k}(x)=&\sum_{j=1}^{K_3}\eta_{kj}M_{k+p}B_j(x), \quad 0\leq\eta_{ij}\leq 1,1\leq k\leq q,\\
    M_i=&\frac{\exp(\delta_i)}{\sum_{k=0}^p{\exp(\delta_k)}}, \quad i=1,\ldots,p+q,\\
    \delta_l\sim&N(0, c_1),\textrm{ for }0\leq l\leq p+q,\\
    \beta_{j}\sim & N(0, c_2)\textrm{ for } 1\leq j\leq K_1,\\
    \theta_{kj}\sim& U(0,1)\textrm{ for }1\leq k\leq p, 1\leq j\leq K_2,\\
    \eta_{kj}\sim& U(0,1)\textrm{ for }1\leq k\leq q, 1\leq j\leq K_3.
\end{align*}
Here $B_{j}$'s are the B-spline basis functions. The parameters $\delta_{j}$'s are unbounded. The verification of support condition \ref{eq:parconditiongarch} for the proposed prior is similar.

\subsection{tviGARCH Model}\label{ssc:tvigarchmodel}
Although the GARCH(1,1) remains one of the most popular models to analyze econometric datasets, empirical evidence shows that these datasets regularly raise suspicion to the parameter space restriction $\sum_i a_i+\sum_j b_j<1$. Note that we used a time-varying analog of this restriction for the tvGARCH modeling in Section \ref{ssc:tvgarchmodel}. This often creates a very stringent condition as the validity of $\sum_i a_i(t)+\sum_j b_j(t)<1$ is questionable. The special case for a time-constant GARCH model where this restriction fails is called an iGARCH model in the literature. We consider the following time-varying analog of iGARCH.

\begin{align}
    &X_i|\mathcal{F}_{i-1}\sim \mathrm{N}(0,\sigma_i^2) \nonumber\\
    &\sigma_i^2=\mu(i/n) + \sum_{k=1}^p a_k(i/n) X_{i-k}^2+\sum_{j=1}^q b_j(i/n) \sigma^2_{i-j}\label{modeligarch}
\end{align}

We impose the following constraints on parameter space for the time-varying parameters, 
\begin{align}\label{eq:parconditionigarch}
    \mathcal{P}=\{\mu, a_i:\mu(x)\geq 0, 0\leq a_{k}(x)\leq 1, \sum_{k}a_{k}(x)+\sum_{j}b_{j}(x)=1\}
\end{align}

\noindent The prior functions that allow us to reformulate the problem keeping it consistent with (\ref{eq:parconditionigarch}) is described below:
\begin{align*}
\mu(x) =&\sum_{j=1}^{K_1}\exp(\beta_j)B_j(x)\label{prior31}\\
     a_{k}(x)=&\sum_{j=1}^{K_2}\theta_{kj}M_{k}B_j(x), \quad 0\leq\theta_{kj}\leq 1,1\leq k\leq p,\\
     b_{i}(x)=&\sum_{j=1}^{K_3}\eta_{kj}M_{k+p}B_j(x), \quad 0\leq\eta_{ij}\leq 1,1\leq i\leq (q-1),\\
     b_{q}(x) = &1-\bigg\{\sum_{k=1}^pa_{k}(x)+\sum_{j=1}^{q-1}b_{j}(x)\bigg\},\\
    M_i=&\frac{\exp(\delta_i)}{\sum_{k=0}^{p+q-1}{\exp(\delta_k)}}, \quad i=1,\ldots,p+q-1,\\
    \delta_l\sim&N(0, c_1),\textrm{ for }0\leq l\leq p+q-1,\\
    \beta_{j}\sim & N(0, c_2)\textrm{ for } 1\leq j\leq K_1,\\
    \theta_{kj}\sim& U(0,1)\textrm{ for }1\leq k\leq p, 1\leq j\leq K_2,\\
    \eta_{kj}\sim& U(0,1)\textrm{ for }1\leq k\leq (q-1), 1\leq j\leq K_3.
\end{align*}

\section{Posterior computation and Implementation}
 \label{sec:comp}
 \subsection{tvARCH structure}
The complete likelihood $L$ of the proposed Bayesian method is given by
\begin{align*}
    L&\propto \exp\bigg(\sum_{i=p}^n \big[-\{\mu(i/n) + \sum_{k=1}^p a_k(i/n) X_{i-k}^2\big\} + X_i^2\log \big\{\mu(i/n) \\
    &\quad+ \sum_{i=1}^p a_i(i/n) X_{i-i}^2\}\big] - \sum_{j=1}^{K_1} \beta_j^2/(2c_2) - \sum_{l=0}^p \delta_l^2/(2c_1)\bigg){\mathbf 1}_{0\leq\theta_{kj}\leq 1},
\end{align*}
where $\mu(x) =\sum_{j=1}^{K_1}\exp(\beta_j)B_j(x), a_{k}(x)=\sum_{j=1}^{K_2}\theta_{kj}M_{k}B_j(x)$ and $M_j=\frac{\exp(\delta_j)}{\sum_{k=0}^p{\exp(\delta_k)}}$. We develop efficient Markov Chain Monte Carlo (MCMC) algorithm to sample the parameter $\beta,\theta$ and $\delta$ from the above likelihood. The computation of derivatives allows us to develop an efficient gradient-based MCMC algorithm to sample these parameters. We calculate the gradients of negative log-likelihood $(-\log L)$ with respect to the parameters $\beta$, $\theta$ and $\delta$. The gradients are given below,
\begin{align}
    &-\frac{d\log L}{\beta_j}=\exp(\beta_j)\bigg(1-\sum_i  \frac{B_j(i/n)X_i^2}{(\mu(i/n)+\sum_{j}a_{j}(i/n)X_{i-j}^2)}\bigg) + \beta_j/c_2,\\
    &-\frac{d\log L}{\theta_{kj}}=M_{k}\bigg(1-\sum_i \frac{B_{j}(i/n)X_i^2}{(\mu(i/n)+\sum_{j}a_{j}(i/n)X_{i-j}^2)}\bigg)&,\\
    &-\frac{d\log L}{\delta_j}=\delta_j/c_1+\sum_k (M_j{\mathbf 1}_{\{j=k\}}-M_jM_k)\sum_i\theta_{kj}B_j(x)\bigg(1-\sum_i \frac{B_{j}(i/n)X_i^2}{(\mu(i/n)+\sum_{j}a_{j}(i/n)X_{i-j}^2)}\bigg),
\end{align}
where ${\mathbf 1}_{\{j=k\}}$ stands for the indicator function which takes the value 1 when $j=k$. 

\subsection{tvGARCH / tviGARCH structure}
The complete likelihood $L_2$ of the proposed Bayesian method of~\eqref{modelgarch} is given by
\begin{align*}
    L_2&\propto \exp\bigg(\sum_{t=p}^n \big[-\{\mu(i/n) + \sum_{i=1}^p a_i(i/n) X_{t-i}+\sum_{i=1}^q b_i(i/n) \lambda_{t-i}\big\} + X_t\log \big\{\mu(i/n) \\
    &\quad+ \sum_{i=1}^p a_i(i/n) X_{t-i}+\sum_{i=1}^q b_i(i/n) \lambda_{t-i}\}\big] - \sum_{j=1}^{K_1} \beta_j^2/(2c_2) - \sum_{l=0}^p \delta_l^2/(2c_1)\\
    &\quad-(d_1+1)\log\lambda_0 -d_1/\lambda_0\bigg){\mathbf 1}_{0\leq\theta_{ij},\eta_{ij}\leq 1},
\end{align*}
We calculate the gradients of negative log-likelihood $(-\log L_2)$ with respect to the parameters $\beta$, $\theta$, $\eta$ and $\delta$. The gradients are given below,
\begin{align*}
    &-\frac{d\log L_2}{\beta_j}=\exp(\beta_j)\bigg(1-\sum_t  \frac{B_j(i/n)X^2_{i-j}}{(\mu(i/n)+\sum_{j}a_{j}(i/n)X^2_{i-j})+\sum_{k}b_{k}(i/n)\sigma^2_{i-k})}\bigg) + \beta_j/c_2,\\
    &-\frac{d\log L_2}{\theta_{lj}}=M_{l}\bigg(1-\sum_t \frac{B_{j}(i/n)X^2_{i-j}}{(\mu(i/n)+\sum_{j}a_{j}(i/n)X^2_{i-j})+\sum_{k}b_{k}(i/n)\sigma^2_{i-k})}\bigg)&,\\
    &-\frac{d\log L_2}{\eta_{kj}}=M_{p+k}\bigg(1-\sum_t \frac{B_{j}(i/n)\sigma^2_{i-j}}{(\mu(i/n)+\sum_{j}a_{j}(i/n)X^2_{i-j})+\sum_{k}b_{k}(i/n)\sigma^2_{i-k})}\bigg)&,\\
    &-\frac{d\log L_2}{\delta_j}=\delta_j/c_1+\sum_k (M_j{\mathbf 1}_{\{j=k\}}-M_jM_k)\times \nonumber\\&\Bigg[\quad\sum_{i\leq p}\theta_{ij}B_j(x)\bigg(1-\sum_t \frac{B_{j}(i/n)X^2_{i-j}}{(\mu(i/n)+\sum_{j}a_{j}(i/n)X^2_{i-j})+\sum_{k}b_{k}(i/n)\sigma^2_{i-k})}\bigg){\mathbf 1}_{\{j\leq p\}}+\nonumber\\
    &\quad\sum_{1\leq k\leq q}\eta_{kj}B_j(x)\bigg(1-\sum_t \frac{B_{j}(i/n)\sigma^2_t}{(\mu(i/n)+\sum_{j}a_{j}(i/n)X^2_{i-j})+\sum_{k}b_{k}(i/n)\sigma^2_{i-k})}\bigg){\mathbf 1}_{\{j > p\}}\Bigg].
\end{align*}
While fitting tvGARCH$(p,q)$, we assume for any $t<0$, $X^2_t=0,\sigma^2_t=0$. Thus, we need to additionally estimate the parameter $\sigma^2_0$. The derivative of the likelihood concerning $\sigma^2_0$ is calculated numerically using the {\tt jacobian} function from R package {\tt pracma}. For the tviGARCH, the derivatives are similar so we avoid computing them for the sake of brevity.

Based on these gradient functions, we develop gradient-based Hamiltonian Monte Carlo (HMC) sampling. Note that, parameter spaces of $\theta_{kj}$'s have bounded support. We circumvent this by mapping any Metropolis candidate falling outside the parameter space back to the nearest boundary. HMC has two parameters, required to be specified. These are the leap-frog step and the step-size parameter. It is difficult to tune both of them simultaneously. We choose to tune the step size parameter to maintain an acceptance range between 0.6 to 0.8. After every 100 iterations, the step-length is adjusted (increased or reduced) accordingly if it falls outside. \cite{neal2011mcmc} showed that a higher leapfrog step is better for estimation accuracy at the expense of greater computation. To maintain a balance between accuracy and computational complexity, we keep it fixed at 30 and obtain good results.

\section{Large-sample properties}\label{sec:theo}
We now focus on obtaining posterior contraction rates for our proposed Bayesian models. The posterior consistency is studied in the asymptotic regime of increasing number of time points $n$.  We study the posterior consistency with respect to the average Hellinger distance on the coefficient functions which is 
$$d_{1,n}^2=\frac{1}{n}d^2_{H}(\kappa_1,\kappa_2)=\frac{1}{n}\int(\sqrt{f_1}-\sqrt{f_2})^2,$$
where $f_1=\prod_{i=1}^nP_{\kappa_1}(X_{i}|X_{i-1})$ and $f_2$ denotes the corresponding likelihoods.

\noindent \textbf{Definition:}
For a sequence $\epsilon_n$ if $\Pi_n(d(f,f_0)|X^{(n)}\geq M_n\epsilon_n|X^{(n)})\rightarrow 0$ in $F_{\kappa_0}^{(n)}$-probability for every sequence $M_n\rightarrow\infty$, then the sequence $\epsilon_n$ is called the \textit{posterior contraction rate.} 

All the proofs are postponed to the supplementary materials. The proof is based on the general contraction rate result for independent non-i.i.d. observations \citep{ghosal2017fundamentals} and some results on B-splines based finite random series. The exponentially consistent tests are constructed leveraging on the famous Neyman-Pearson Lemma as in \cite{ning2020Bayesian}. Thus the first step is to calculate posterior contraction rate with respect to average log-affinity $r_n^2(f_1,f_2)=-\frac{1}{n}\log\int f_1^{1/2}f_2^{1/2}$. Then we show that $r_n^2(f_1,f_2)\lesssim\epsilon_n^2$ implies $\frac{1}{n}d_H^2(f_1,f_2) \lesssim \epsilon_n^2$. We also consider following simplified priors for $\alpha_j$ and $\tau_i$ to get better control over tail probabilities,
\begin{align}
    \alpha_j\sim\mathrm{Gamma}(g_1,g_1),\quad \tau_i \sim U(0,1). \label{prioringar:choice1}
\end{align}

\subsection{tvARCH model}\label{sec:tvarch}
Let $\kappa=(\mu, a_1)$ stands for the complete set of parameters. For sake of generality of the method, we put a prior on $K_1$ and $K_2$ with probability mass function given by, 
\begin{eqnarray}\label{eq:bij}
\Pi(K_i=k)=b_{i1}\exp[-b_{i2} k (\log k)^{b_{i3}}],
\end{eqnarray}
for $i=1,2$. These priors have not been considered while fitting the model as it would require computationally expensive reversible jump MCMC strategy.  
The contraction rate will depend on the smoothness of true coefficient functions $\mu$ and $a$ and the parameters $b_{13}$ and $b_{23}$ from the prior distributions of $K_1$ and $K_2$. Let $\kappa_0=(\mu_0, a_{01})$ be the truth of $\kappa$.

\noindent Assumptions(A): There exists constants $M_X>1, 0<M_{\mu}<M_X$ such that,
\begin{itemize}
    \item[(A.1)] The coefficient functions satisfy $\sup _x \mu_0(x)<M_{\mu}$ and $\sup_x a_{01}(x) <1-M_{\mu}/M_{X}$.
    \item[(A.2)] $\inf_x \min(\mu_0(x),a_{01}(x))>\rho$ for some small $\rho>0$.
    \item[(A.3)] $E(X_0^2)<M_X$.
\end{itemize}

\noindent Assumptions (A.1)-(A.3) ensure $$\IE_{\kappa_0}(X_i^2)=\IE_{\kappa_0}(\IE_{\kappa_0}(X_i^2|X_{i-1}))<M_{\mu}+\left(1-\frac{M_{\mu}}{M_X}\right)M_X<M_X$$
by recursion. 

\ignore{
\begin{remark}
For the tvARCH(p) case, (A.1) would read the coefficient functions satisfy $$\sup _x \mu_0(x)<M_{\mu}\text{ and }\sup_x \sum_{i=1}^p a_{0i}(x) <1-M_{\mu}/M_{X}.$$
\end{remark}
}


\begin{theorem}
\label{thm:arch}
Under assumptions (A.1)-(A.3), let the true functions $\mu_0(\cdot)$ and $a_{10}(\cdot)$ be H\"older smooth functions with regularity level $\iota_1$ and $\iota_2$ respectively, then the posterior contraction rate with respect to the distance $d_{1,n}^2$ is $$\max\bigg\{n^{-\iota_1/(2\iota_1+1)} (\log n)^{\iota_1/(2\iota_1+1)+(1-b_{13})/2},n^{-\iota_2/(2\iota_2+1)} (\log n)^{\iota_2/(2\iota_2+1)+(1-b_{23})/2}\bigg\},$$
\end{theorem}
\noindent where $b_{ij}$ are specified in (\ref{eq:bij}). 

\subsection{tvGARCH model}\label{sec:tvgarch}
Let $\kappa=(\mu, a_1, b_1)$ stands for the complete set of parameters. For sake of generality of the method, we put a prior on $K_1$, $K_2$ and $K_3$ with probability mass function given by, 
\begin{eqnarray}\label{eq:bij2}
\Pi(K_i=k)=b_{i1}\exp[-b_{i2} k (\log k)^{b_{i3}}],
\end{eqnarray}
for $i=1,2$. These priors have not been considered while fitting the model as it would require computationally expensive reversible jump MCMC strategy.  
The contraction rate will depend on the smoothness of true coefficient functions $\mu$ and $a$ and the parameters $b_{13}$ and $b_{23}$ from the prior distributions of $K_1$ and $K_2$. Let $\kappa_0=(\mu_0, a_{01})$ be the truth of $\kappa$.

\noindent Assumptions(B): There exists constants $M_X>1, 0<M_{\mu}<M_X$ such that,
\begin{itemize}
    \item[(B.1)] The coefficient functions satisfy $\sup _x \mu_0(x)<M_{\mu}$ and $\sup_x (a_{01}(x)+b_{01}(x)) <1-M_{\mu}/M_{X}$.
    \item[(B.2)] $\inf_x \min(\mu_0(x),a_{01}(x), b_{01}(x))>\rho$ for some small $\rho>0$.
    \item[(B.3)] $E(X_0^2)<M_X, \sigma_{00}^2<M_X$. 
\end{itemize}

\noindent Assumptions (B.1) and (B.3) ensure $$\IE_{\kappa_0}(X_i^2)=\IE_{\kappa_0}(\IE_{\kappa_0}(X_i^2|X_{i-1}))<M_{\mu}+\left(1-\frac{M_{\mu}}{M_X}\right)M_X<M_X$$
by recursion. Similarly we have $\IE(\sigma_i^2)=\IE_{\kappa_0}(\IE_{\kappa_0}(X_i^2|\mathcal{F}_i))=\IE_{\kappa_0}(X_i^2)<M_X$.


\begin{theorem}
\label{thm:garch}
Under assumptions (B.1)-(B.3), let the true functions $\mu_0(\cdot)$, $a_{10}(\cdot)$ and $b_{10}(\cdot)$ be H\"older smooth functions with regularity level $\iota_1$. $\iota_2$ and $\iota_3$ respectively, then the posterior contraction rate with respect to the distance $d_{1,n}^2$ is 
\begin{align*}
    \max\bigg\{&n^{-\iota_1/(2\iota_1+1)} (\log n)^{\iota_1/(2\iota_1+1)+(1-b_{13})/2},n^{-\iota_2/(2\iota_2+1)} (\log n)^{\iota_2/(2\iota_2+1)+(1-b_{23})/2},\\&n^{-\iota_3/(2\iota_3+1)} (\log n)^{\iota_2/(2\iota_3+1)+(1-b_{33})/2}\bigg\},
\end{align*}
\end{theorem}
\noindent where $b_{ij}$ are specified in (\ref{eq:bij2}). 

\ignore{
\begin{remark}
For the tvGARCH(p,q) case, (B.1) would read the coefficient functions satisfy $$\sup _x \mu_0(x)<M_{\mu}\text{ and }\sup_x (\sum_{i=1}^p a_{0i}(x)+\sum_{j=1}^q b_{0j})(x) <1-M_{\mu}/M_{X}.$$
\end{remark}
}
\subsection{tviGARCH model}\label{sec:tvigarch}

Let $\kappa=(\mu, a_1)$ stands for the complete set of parameters. For sake of generality of the method, we put a prior on $K_1$ and $K_2$ with probability mass function given by, 
\begin{eqnarray}\label{eq:bij2garch}
\Pi(K_i=k)=b_{i1}\exp[-b_{i2} k (\log k)^{b_{i3}}],
\end{eqnarray}
for $i=1,2$. These priors have not been considered while fitting the model as it would require computationally expensive reversible jump MCMC strategy.  
The contraction rate will depend on the smoothness of true coefficient functions $\mu$ and $a$ and the parameters $b_{13}$ and $b_{23}$ from the prior distributions of $K_1$ and $K_2$. Let $\kappa_0=(\mu_0, a_{01})$ be the truth of $\kappa$.

\noindent 
\begin{itemize}
    \item[(C.1)] The coefficient functions satisfy $\sup _x \mu_0(x)<M_{\mu}<\infty$ for some $M_{\mu}$
    \item[(C.2)] $\inf_x (\mu_0(x))>\rho ,\inf_x a_{01}(x)>\rho, \sup_x a_{0,1}(x)<1-\rho$  for some $\rho>0$.
\end{itemize}


\begin{theorem}
\label{thm:igarch}
Under assumptions (C.1)-(C.2), let the true functions $\mu_0(\cdot)$ and $a_{10}(\cdot)$ be H\"older smooth functions with regularity level $\iota_1$. and $\iota_2$respectively, then the posterior contraction rate with respect to the distance $d_{1,n}^2$ is $$\max\bigg\{n^{-\iota_1/(2\iota_1+1)} (\log n)^{\iota_1/(2\iota_1+1)+(1-b_{13})/2},n^{-\iota_2/(2\iota_2+1)} (\log n)^{\iota_2/(2\iota_2+1)+(1-b_{23})/2}\bigg\},$$
\end{theorem}
\noindent where $b_{ij}$ are specified in (\ref{eq:bij2garch}). 

\section{Simulation}
\label{sec:sim}
We run simulations to study the performance of our proposed Bayesian method in capturing the true coefficient functions under different true models. The hyperparameters $c_1$ and $c_2$ of the normal prior are all set 100, which makes the prior weakly informative. We consider 4, 5 and 6 equidistant knots for the B-splines when $n=200, 500$ and $1000$ respectively. We collect 10000 MCMC samples and consider the last 5000 as post burn-in samples for inferences. We shall compare the estimated functions with the true functions in terms of the posterior estimates of functions along with its 95\% pointwise credible bands. The credible bands are calculated from the MCMC samples at each point $t=1/T,2/T,\ldots, 1$. We take the posterior mean as the posterior estimate of the unknown functions.

Since, to the best of our knowledge, there is no other Bayesian model for these time-varying conditional heteroscedastic models, we compare our Bayesian estimates with corresponding frequentist time-varying estimates. For computing the time-varying estimates of these models, we use the kernel-based method from \cite{sayar2020}. The M-estimator of the parameter vector $\theta(t)$ are obtained using the conditional quasi log-likelihood. For instance, in the tvARCH(1) case, say $\theta(t)=(\mu(t),a_1(t))$
\begin{equation*}
    \hat \theta_{b_n}(t)= \argmin_{\theta \in \Theta} \sum_{i=2}^{n}K\left(\frac{t-i/n}{b_n}\right) \ell(X_i|\mathcal{F}_{i-1},\theta) \quad\quad t \in [0,1].\label{eq:estimator}
\end{equation*}

\noindent where $\ell(\cdot)$ denotes the Gaussian log-likelihood. Note that these methods are fast but usually need a cross-validated choice of bandwidth $b_n$. We use $K(x)=3/4(1-x^2)\mathbf{I}(|x| \leq 1)$ and an appropriately chosen bandwidth as suggested by the authors therein. Since our discussion also involves iGARCH formulation, we wrote a separate kernel-based frequentist estimation for iGARCH models analogously. Apart from these two time-varying estimates, we also obtain a time-constant fit on the same data to help initiate a discussion on whether there was a necessity of introducing coefficients varying with time. For this \texttt{tseries} and \texttt{rugarch} R packages are used respectively for ARCH/GARCH and iGARCH fits. 

To compare these estimates, we evaluate the average mean square errors (AMSE) for the three estimates.
Note that in an usual linear regression of response $y$ on predictor $X$ scenario, the fitted MSE is often defined as $\frac{1}{n}\sum (y_i - \hat{y}_i)^2$. Since here, $X_i|\mathcal{F}_{i-1} \sim N(0,\sigma_i^2)$, we use the following definition of AMSE $$\text{AMSE}=\frac{1}{n}\sum_i(X_i^2-\hat{\sigma}_i^2)^2$$ where the $\hat{\sigma}_i^2$ is computed with the fitted parameter values as per the model under consideration. For example, for a tvGARCH(1,1) model we have  $$\hat{\sigma}_i^2= \hat{\mu}(i/n)+\hat{a}(i/n)X_{i-1}^2+\hat{b}(i/n)\hat{\sigma}_{i-1}^2$$ where $\hat{\mu}(\cdot), \hat{a}(\cdot)$ and $\hat{b}(\cdot)$ are the estimated curves from the posterior. Replacing the response $y_i$ by $X_i^2$ is natural as often autocorrelations of $X_i^2$ are checked to gauge presence of CH effect. Moreover, one of the early methods to deal with CH models was to view $X_i^2$ approximated by an TVAR(1) process. See \cite{bose2003estimating} and references therein. Similar estimators as our proposed AMSE to evaluate the fitting accuracy has been used in the literature previously. See \cite{starica2003garch, fry08, tvgarch2013, sayar2020} for example. 

In the next three subsections, we provide the results for the three models, namely, tvARCH, tvGARCH, and tviGARCH. Our conclusions from these results are illustrated at the end of the section.

\subsection{tvARCH case}

We start by considering the following tvARCH(1) model from \ref{modelarch}. Three different choices for $n$ are considered, $n=200,500$ and $1000$. The true functions are,
\begin{align*}
    \mu_0(x)=&10\exp\big(-(x-0.5)^2/0.1\big),\\
    a_{10}(x)=&0.4(x-0.15)^2+0.1.
\end{align*}

\begin{figure}[htbp]
\centering
        \subfigure[$n=200$]{\label{fig:c.11}\includegraphics[width=130mm, height=30 mm]{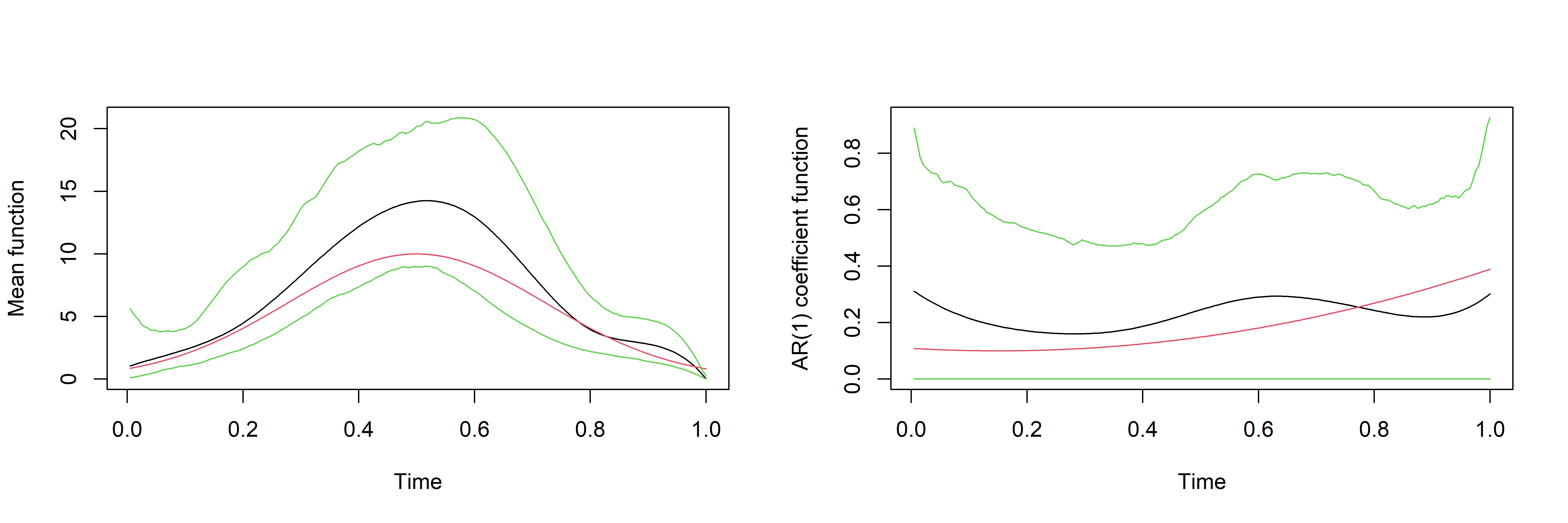}}
        \subfigure[$n=500$]{\label{fig:c.12}\includegraphics[width=130mm, height=30 mm]{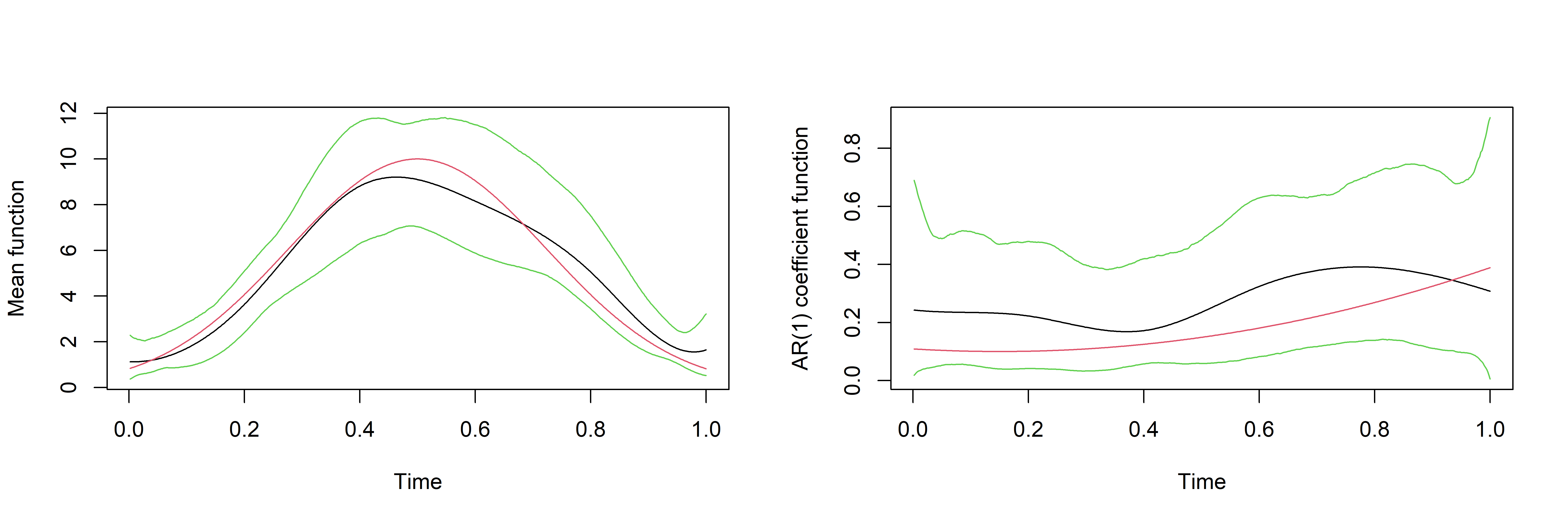}}
        \subfigure[$n=1000$]{\label{fig:c.13}\includegraphics[width=130mm, height=30 mm]{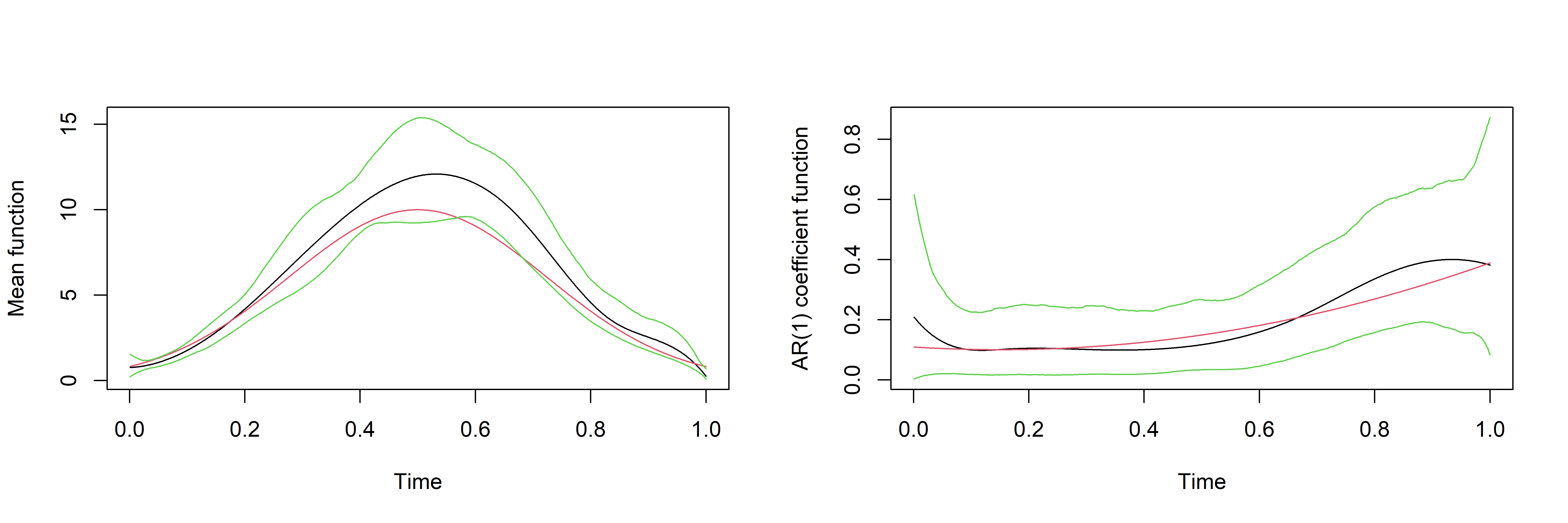}}
        \caption{tvARCH(1):True coefficient functions(red), estimated curve(black) along with the 95\% pointwise credible bands (green) are shown for T=200,500,1000 from top to bottom} 
        \label{archpic}
\end{figure}

\noindent We compare the estimated functions with the truth for sample size 1000 in Figures~\ref{archpic}. Table~\ref{AMSEtvARCH1} illustrates the performance of our method with respect to other competing methods.

\begin{table}[ht]
\centering
\caption{AMSE comparison for different sample sizes across different methods when the true model is tvARCH with $p=1$}
\begin{tabular}{rrrr}
  \hline
 & ARCH(1) &  Frequentist tvARCH(1) & Bayesian tvARCH(1) \\ 
  \hline
   $n=200$ & 96.42 & 90.34 & \textbf{85.22} \\ 
  $n=500$ & 128.07 & 122.53 & \textbf{118.45} \\ 
  $n=1000$ & 138.06 & 130.33 & \textbf{127.06} \\ 
   \hline
\end{tabular}
\label{AMSEtvARCH1}
\end{table}

\subsection{tvGARCH case}
Next we explore the following GARCH(1,1) model (cf. \ref{modelgarch})for different choices of $n$. The true functions are,
\begin{align*}
    \mu_0(x)=&1-0.8\sin(\pi x/2),\\
    a_{10}(x)=&0.5-(x-0.3)^2,\\
    b_{10}(x)=&0.4-0.5(x-0.4)^2
\end{align*}

\begin{figure}[htbp]
\centering
        \subfigure[$n=200$]{\label{fig:c.21}\includegraphics[width=130mm,height=30mm]{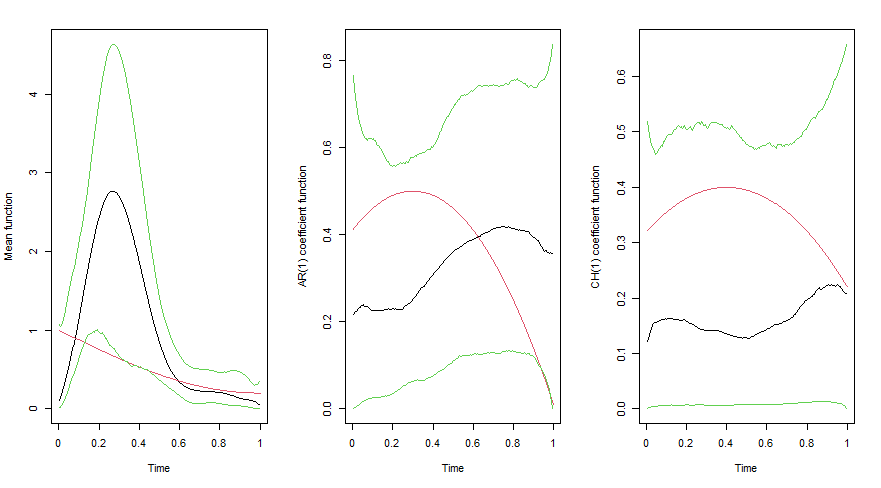}}
        \subfigure[$n=500$]{\label{fig:c.22}\includegraphics[width=130mm,height=30mm]{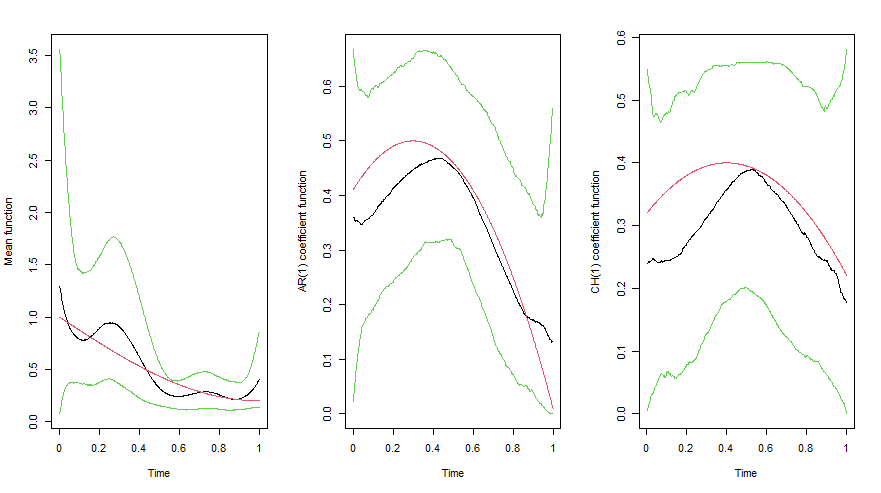}}
        \subfigure[$n=1000$]{\label{fig:c.23}\includegraphics[width=130mm,height=30mm]{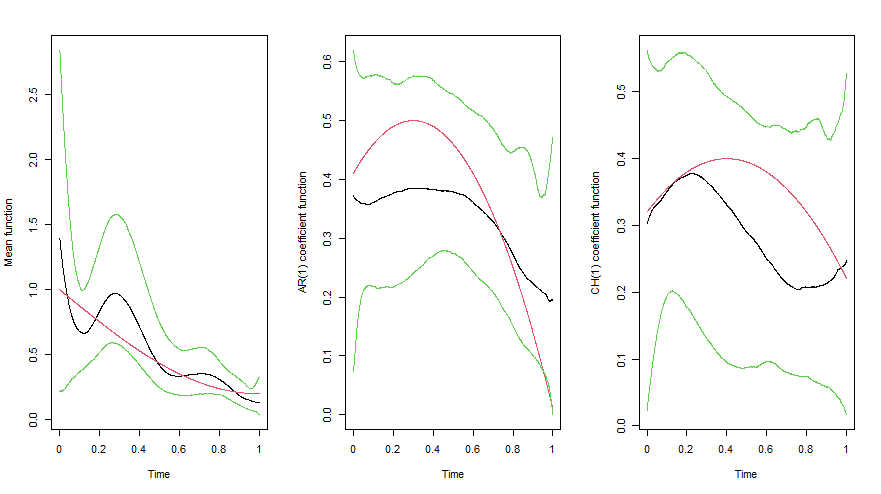}}
        \caption{tvGARCH(1,1):True coefficient functions(red), estimated curve(black) along with the 95\% pointwise credible bands (green) are shown for T=200,500,1000 from top to bottom} 
        \label{garchpic}
\end{figure}

Note that, estimation of GARCH, due to the additional $b_i(\cdot)$ parameter curves is a significantly more challenging problem and often requires a much higher sample size to have a reasonable estimation. We show by the means of the following pictures in Figure \ref{garchpic} that the estimation looks reasonable even for smaller sample sizes. The AMSE score comparisons are shown in Tables~\ref{AMSEtvGARCH1}. The performance of our method is also contrasted with other competing methods.

\begin{table}[ht]
\centering
\caption{AMSE comparison for different sample sizes across different methods when the true model is tvGARCH(1,1)}
\begin{tabular}{rrrr}
  \hline
 & GARCH(1,1) &  Frequentist tvGARCH(1,1) & Bayesian tvGARCH(1,1) \\ 
  \hline
    $n=200$ & 33.99 & 31.84 & \textbf{29.43} \\  
$n=500$ & 45.46 & 34.77  & \textbf{33.33} \\ 
 $n=1000$ & 42.60 & 37.09 &  \textbf{36.55} \\ 
   \hline
\end{tabular}
\label{AMSEtvGARCH1}
\end{table}

\subsection{tviGARCH case}
Finally we consider the tviGARCH(1,1) model (cf. \ref{modeligarch}) a special case of GARCH. Note that due to the constraint $a_1(\cdot)+b_1(\cdot)=1$ we only consider the mean function and AR(1) function for plotting. For this case, our true functions are as follows 
\begin{align*}
    \mu_0(x)=&\exp\big(-(x-0.5)^2/0.1\big),\\
    a_{10}(x)=&0.4(x-1)^2+0.1.
\end{align*}

\begin{figure}[htbp]
\centering
        \subfigure[$n=200$]{\label{fig:c.31}\includegraphics[width=130mm, height=30 mm]{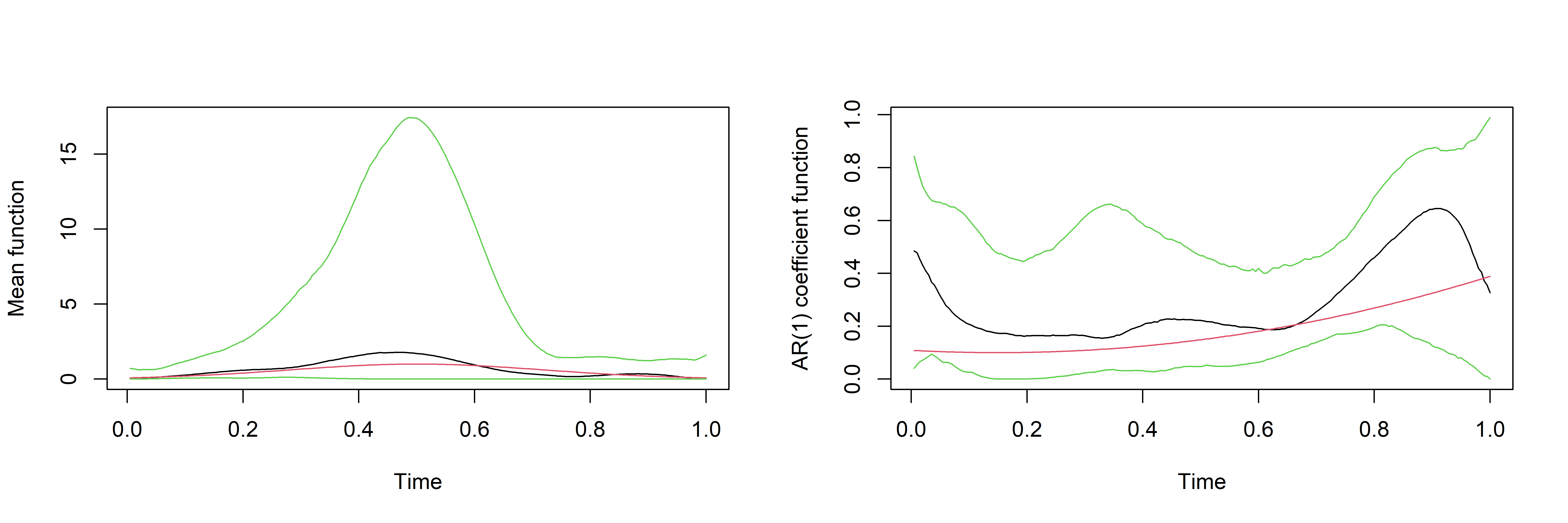}}
        \subfigure[$n=500$]{\label{fig:c.32}\includegraphics[width=130mm, height=30 mm]{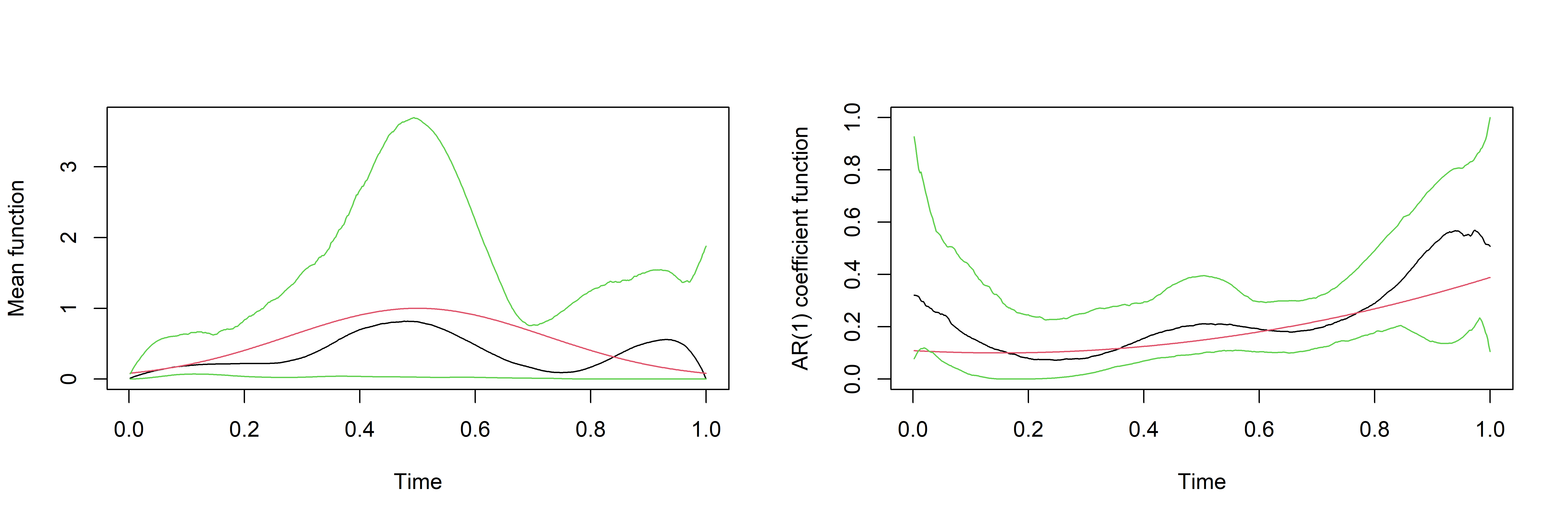}}
        \subfigure[$n=1000$]{\label{fig:c.33}\includegraphics[width=130mm, height=30 mm]{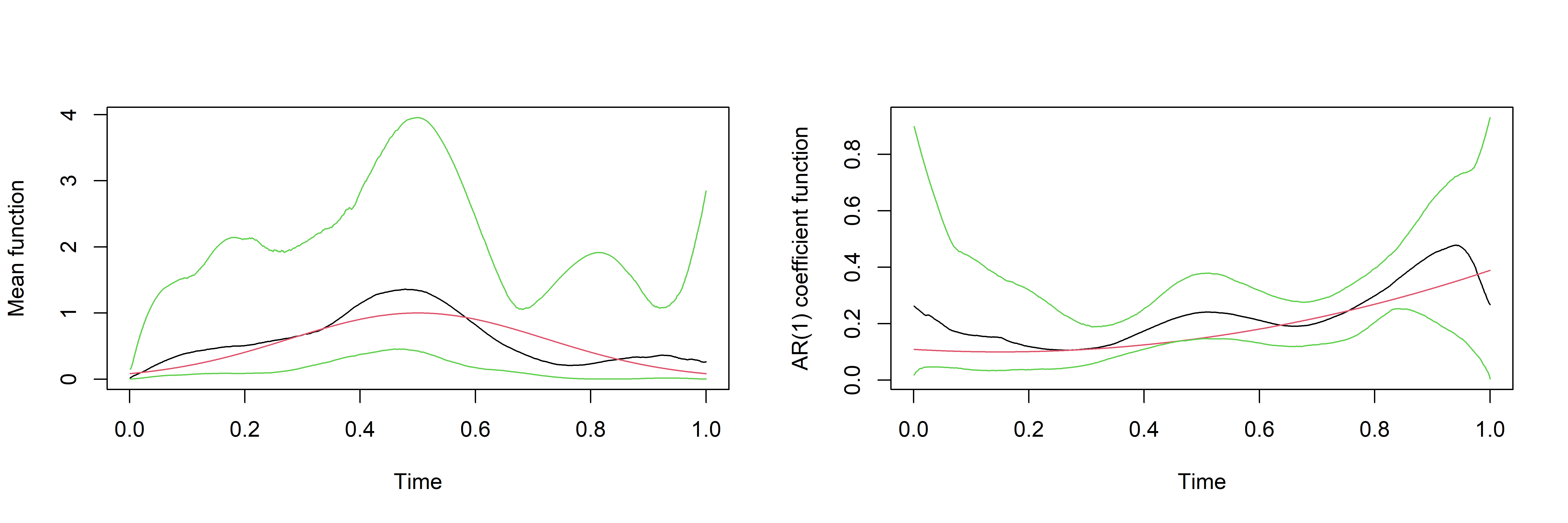}}
        \caption{tviGARCH(1,1):True coefficient functions(red), estimated curve(black) along with the 95\% pointwise credible bands (green) are shown for T=200,500,1000 from top to bottom} 
        \label{igarchpic}
\end{figure}

\noindent  The frequentist computation for tviGARCH method is carried out based on a kernel-based estimation scheme along the same line as \cite{sayar2020}. The estimated plots along with the 95\% credible intervals are shown in Figure \ref{igarchpic} for three sample sizes $n=200,500,1000$ and the AMSE scores in Table~\ref{AMSEtviGARCH1}. 

\begin{table}[ht]
\centering
\caption{AMSE* comparison for different sample sizes across different methods when the true model is tviGARCH with $p=1,q=1$. AMSE* stands for mean of the log(AMSE)}
\begin{tabular}{rrrr}
  \hline
 & iGARCH(1,1) &  Frequentist tviGARCH(1,1) & Bayesian tviGARCH(1,1) \\ 
  \hline
  200 & 8.20 & 23.86 & \textbf{8.14} \\ 
  500 & 9.06 & 18.72 & \textbf{9.06} \\ 
  1000 & 10.59 & 25.92 & \textbf{10.59} \\
   \hline
\end{tabular}
\label{AMSEtviGARCH1}
\end{table}

 To summarize, our estimated functions are close to true functions for all the cases. We also find that the credible bands are getting tighter with increasing sample size. Thus estimates are improving in precision as sample size increases as shown in Figures~\ref{archpic} to~\ref{igarchpic}. AMSEs of our Bayesian estimates are at least better for all the cases as in Tables~\ref{AMSEtvARCH1} to~\ref{AMSEtviGARCH1}. For tviGARCH, AMSE* is considered due to the huge and somewhat incomparable values of AMSE due to non-existent variance.

\ignore{
 $n=200$& 18472.35 & 11979794899825756417884628.00 & 19204.45 \\ 
  $n=500$ & 388817.06 & 1355052130001898700800.00 & 385843.59 \\ 
  $n=1000$ & 500533.23 & 911317428000025424242488822.00 & 488315.33 \\

}

\section{Real data application}
\label{sec:application}
Towards applying our methods on real-life datasets we stick to econometric data for varying time horizons. These datasets show considerable time-variation justifying our models to be suitable for understanding how the parameter functions have evolved.
Typically we model the log-return data of the daily closing price of these data to avoid the unit-root scenario. The log-return is defined as follows and is close to the relative return

$$Y_i=\log P_i- \log P_{i-1}=\log \left( 1+ \frac{P_i-P_{i-1}}{P_{i-1}}\right) \approx \frac{P_i-P_{i-1}}{P_{i-1}}, $$
where $P_i$ is the closing price on the $i^{th}$ day. Conditional heteroscedastic models are popularly used for model building, analysis and forecasting. Here we extend such models to a more sophisticated and general scenario by allowing the coefficient functions to vary. 

In this section, we analyze two datasets: USD to JPY conversion and NASDAQ, a popular US stock market data. We analyze the NASDAQ data through tvGARCH(1,1) and tviGARCH(1,1) models and USDJPY conversion rate data through tvARCH(1) models. We just fit one lag for these models as multiple lag fits are similar and larger lags seem to be insignificant. This result is consistent with the findings in \cite{sayar2020}, \cite{subbarao2008} etc. Moreover, as \cite{subbarao2008} claims, stock indices and Forex rates are more suited to GARCH and ARCH type of models respectively for their superior predictive performance.  Each of these datasets was collected up to 31 July 2020. We exhibit our results for the last 200,500 and 1000 days which capture the last 6 months, around 1.5 years, and around 3 years of data respectively. All these datasets were collected from \url{www.investing.com}. Note that these datasets are usually available for weekdays barring holidays and typically there are about 260 data points every single year.

\subsection{USDJPY data: tvARCH(1) model}\label{ssc:usjpy}

We obtain the following Figure \ref{fig:archdata} that shows our estimation for fitting a tvARCH(1) model on the USD to JPY conversion data for the last 200,500 and 1000 days ending on 31 July 2020. The AMSE is also computed and contrasted with other competing methods in Table \ref{AMSEtvARCH1data}.  Figure~\ref{fig:archdata} depicts the estimated functions with 95\% credible bands for different sample sizes. One can see the bands become much shorter for larger sample sizes. The mean coefficient function $\mu(\cdot)$ is generally time-varying for all three cases as one cannot fit a horizontal line through the 95\% posterior bands. There seems to be a rise in the mean value around 100 days ago from July 31, 2020, which is right around the time the COVID-19 pandemic hit the world. With the analysis of $n=1000$ days, we see that the volatility is quite high around October 2016 which coincides with the presidential election time of 2016. The AR(1) coefficient does not show the huge time-varying property. We also tabulate the AMSE for the three sample sizes in Table~\ref{AMSEtvARCH1data} and one can see for smaller sample sizes such as $n=200$, the proposed Bayesian tvARCH achieves a significantly lower score but when the sample size grows then the performance becomes similar.

\begin{figure}[htbp]
\centering
        \subfigure[$n=200$]{\label{fig:c.41}\includegraphics[width=130mm, height=40 mm]{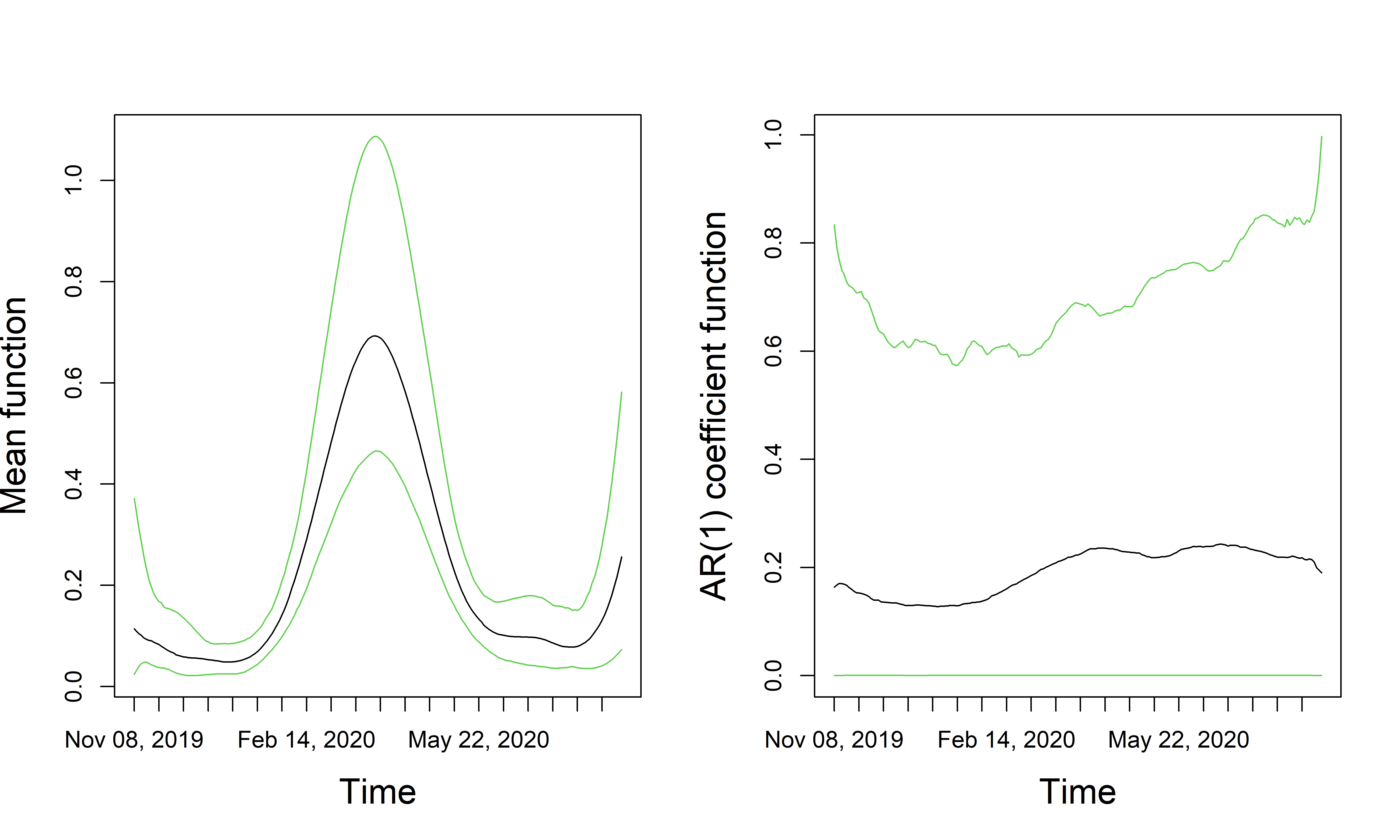}}
        \subfigure[$n=500$]{\label{fig:c.42}\includegraphics[width=130mm, height=40 mm]{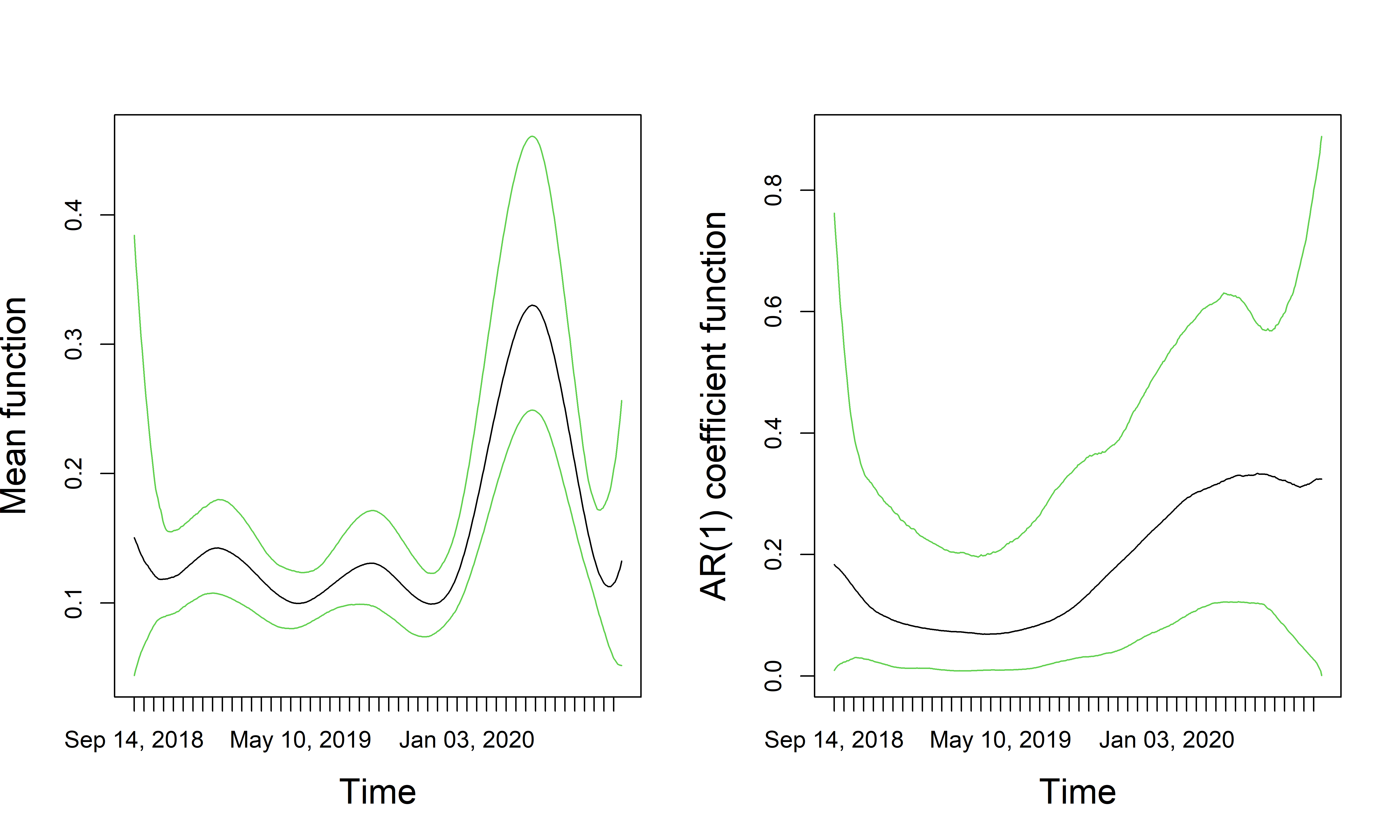}}
        \subfigure[$n=1000$]{\label{fig:c.43}\includegraphics[width=130mm, height=40 mm]{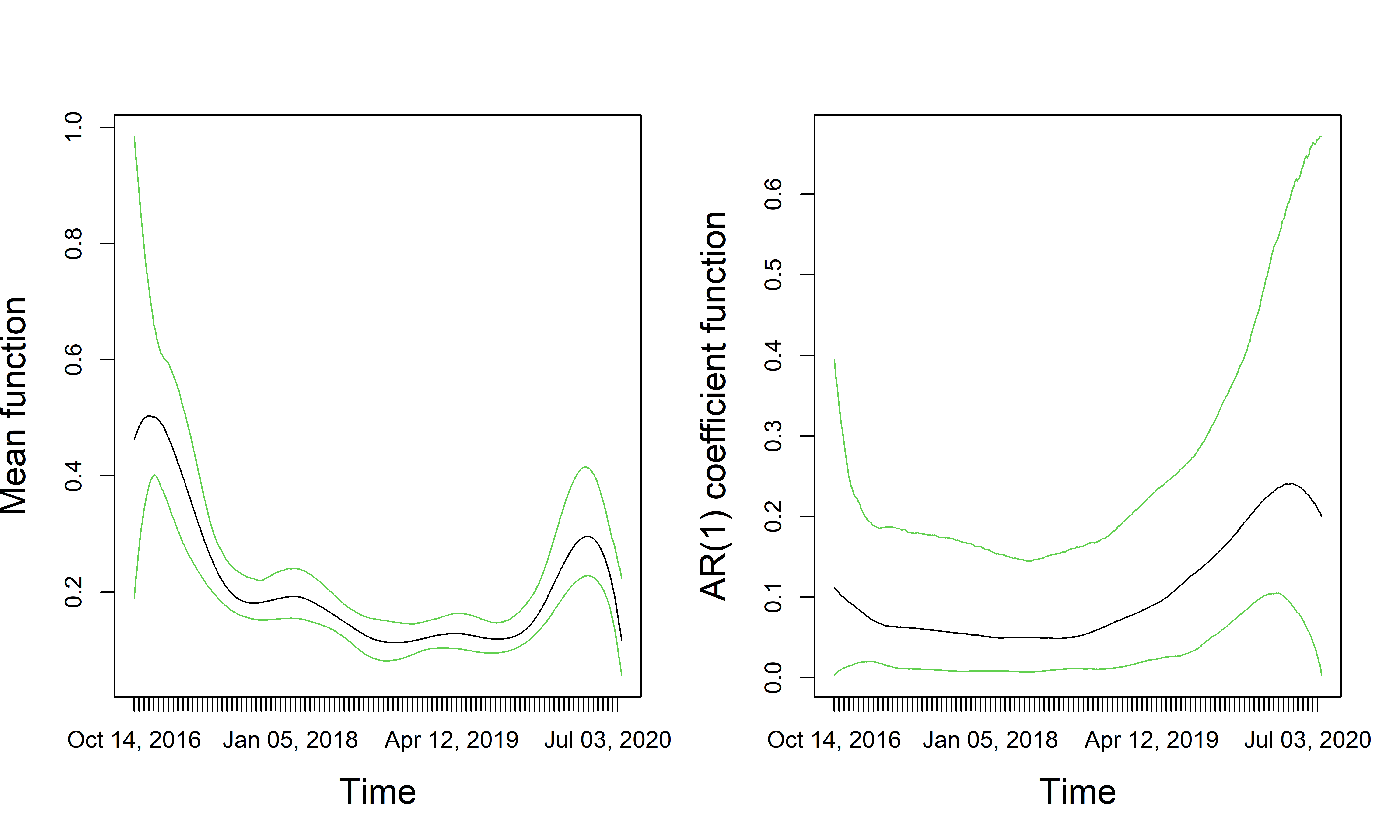}}
        \caption{USDJPY data (tvARCH(1) model) Estimated curve(black) along with the 95\% pointwise credible bands (green) are shown for T=200,500,1000 from top to bottom} 
        \label{fig:archdata}
\end{figure}

\begin{table}[ht]
\centering
\caption{AMSE comparison: tvARCH(1) model- USDJPY data}
\begin{tabular}{rrrr}
  \hline
 & ARCH(1) &  Frequentist tvARCH(1) & Bayesian tvARCH(1) \\ 
  \hline
    $n=200$ & 1.4572 & 1.2259 & \textbf{1.1712} \\  
$n=500$ & 0.6281 & 0.5313 & \textbf{0.5218} \\ 
 $n=1000$ & 0.4265 & \textbf{0.3773} & 0.3785 \\ 
   \hline
\end{tabular}
\label{AMSEtvARCH1data}
\end{table}

\subsection{NASDAQ data: tvGARCH(1,1) model}\label{ssc:nasdaqgarch}

As has become standard in analyzing stock market datasets using GARCH models, we use time-varying GARCH for small orders. We obtain the following Figure \ref{fig:garchdata} for fitting a tvGARCH(1,1) model on the NASDAQ data for the last 200,500, 1000 days ending on 31 July 2020. One can see the $a_1(\cdot)$ values are generally low and the $b_1(\cdot)$ values are higher which is consistent with how these outcomes turn out for time-constant estimates for econometric datasets. One can also see the role sample size plays in curating these time-varying estimates. For $n=200$, the $b_1(\cdot)$ achieves high value of 0.6 around mid-March 2020 but for higher sample sizes it shows values as high as 0.8. One can also note the striking similarity for the analysis of the last 500 and 1000 days which is fairly consistent with the idea that estimation is more stable for such CH type models with a higher sample size. Nonetheless, the estimates for $n=200$ seem quite smooth as well which can be seen as a benefit of our methodology.
Table~\ref{AMSEtvGARCH1data} provides a comparison of AMSE scores across the three methods for three sample sizes. The Bayesian tvGARCH(1,1) performs relatively better than other methods and estimated curves have smaller credible bands with a growing sample size. The behavior of the mean function also shows higher volatility around the pandemic.

\begin{figure}[htbp]
\centering
        \subfigure[$n=200$]{\label{fig:c.51}\includegraphics[width=130mm, height=40 mm]{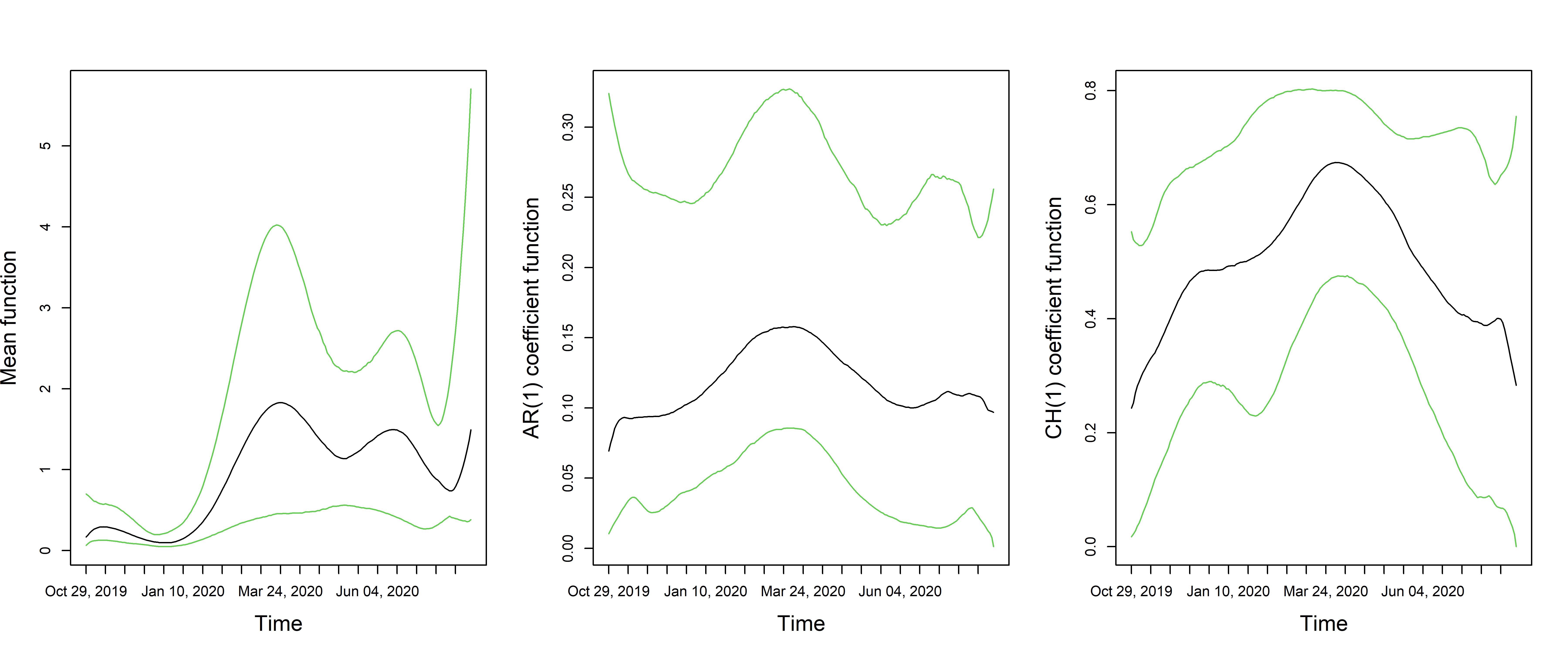}}
        \subfigure[$n=500$]{\label{fig:c.52}\includegraphics[width=130mm, height=40 mm]{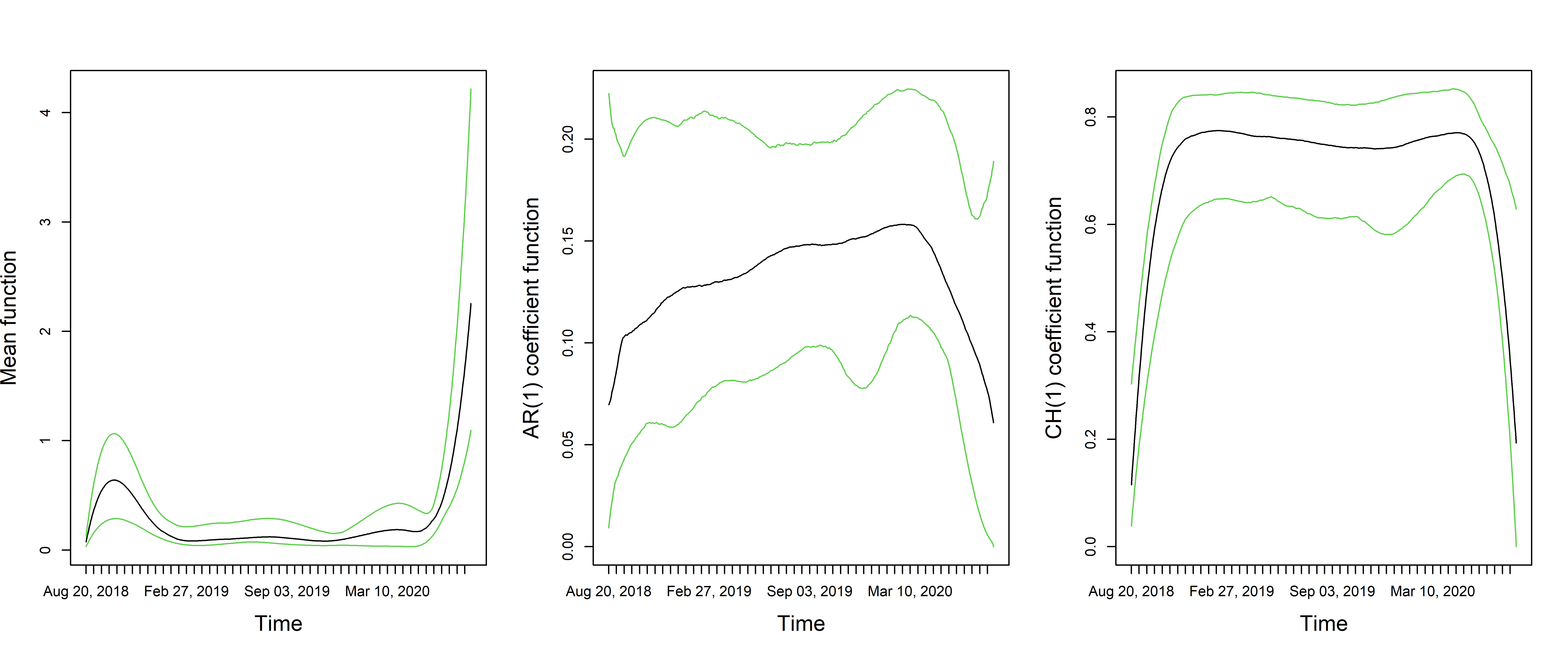}}
        \subfigure[$n=1000$]{\label{fig:c.53}\includegraphics[width=130mm, height=40 mm]{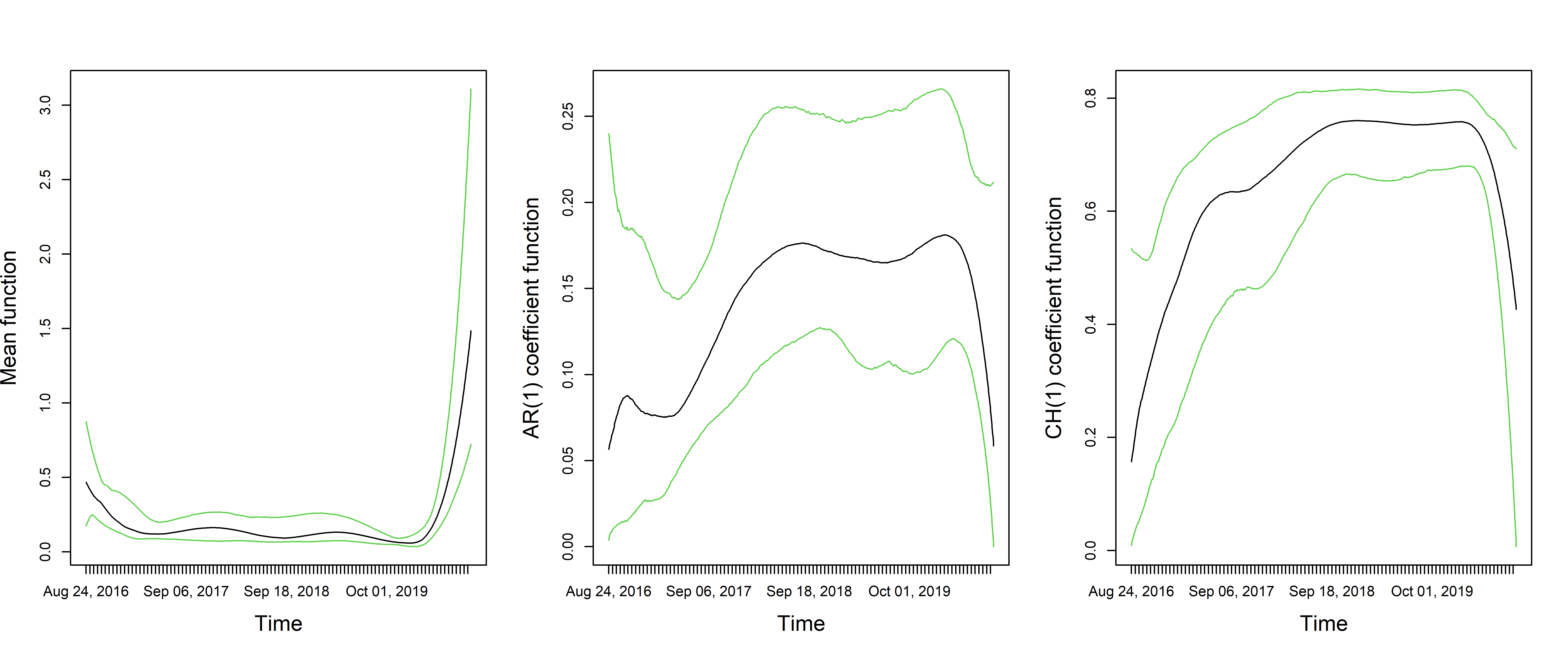}}
        \caption{NASDAQ data (tvGARCH(1,1) model) Estimated curve(black) along with the 95\% pointwise credible bands (green) are shown for T=200,500,1000 from top to bottom} 
        \label{fig:garchdata}
\end{figure}

\begin{table}[ht]
\centering
\caption{AMSE comparison: tvGARCH(1,1) model- NASDAQ data}
\begin{tabular}{rrrr}
  \hline
 & GARCH(1,1) &  Frequentist tvGARCH(1,1) & Bayesian tvGARCH(1,1) \\ 
  \hline
    $n=200$ & 203.5917 & 203.5917 & \textbf{202.6192} \\
$n=500$ & 104.7443 & 90.5395 & \textbf{90.3126} \\ 
 $n=1000$ & 46.16759 & 46.9225 & \textbf{45.5618} \\ 
   \hline
\end{tabular}
\label{AMSEtvGARCH1data}
\end{table}

\subsection{NASDAQ data: tviGARCH(1,1) model}\label{ssc:nasdaqigarch}
\begin{figure}[htbp]
\centering
        \subfigure[$n=200$]{\label{fig:c.61}\includegraphics[width=100mm, height=40 mm]{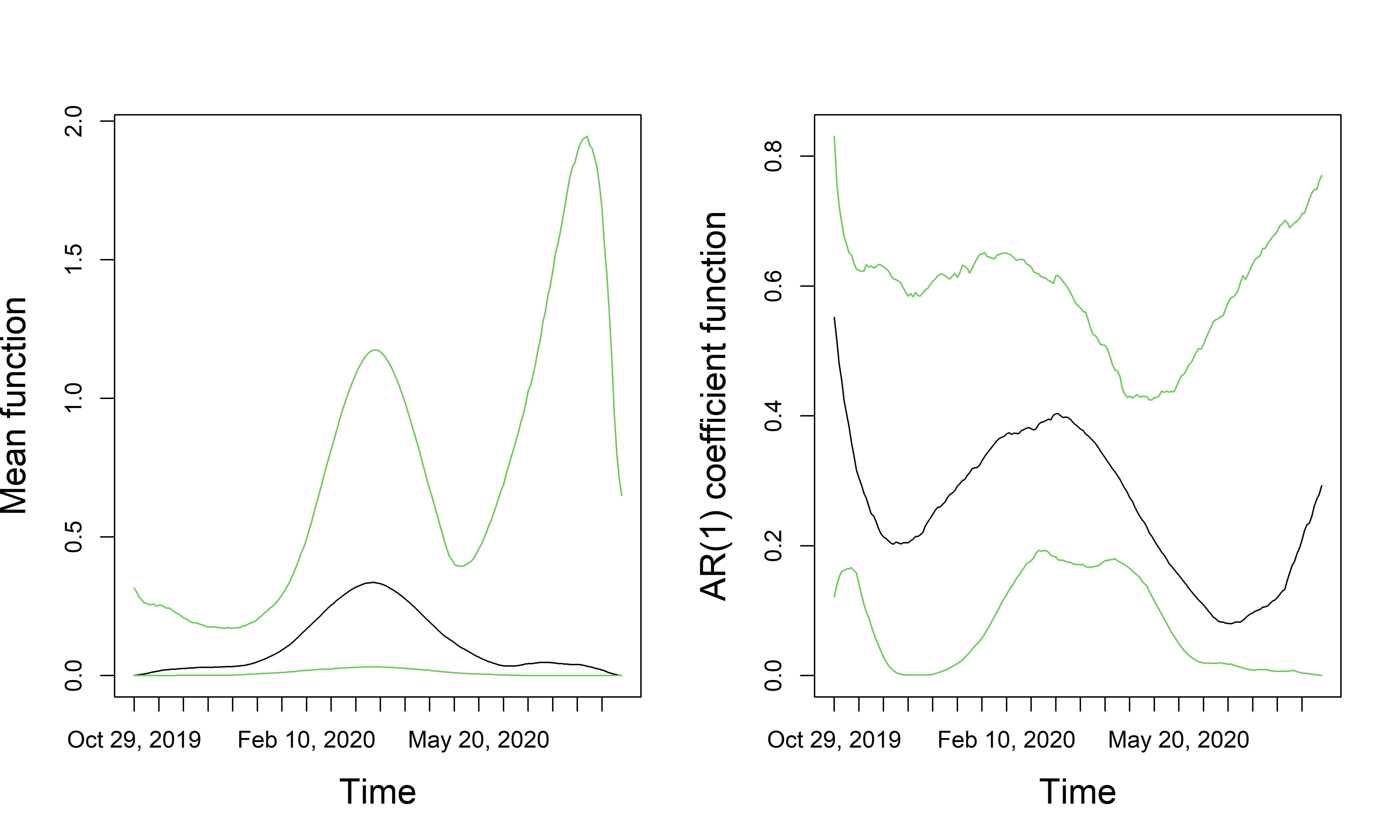}}
        \subfigure[$n=500$]{\label{fig:c.62}\includegraphics[width=100mm, height=40 mm]{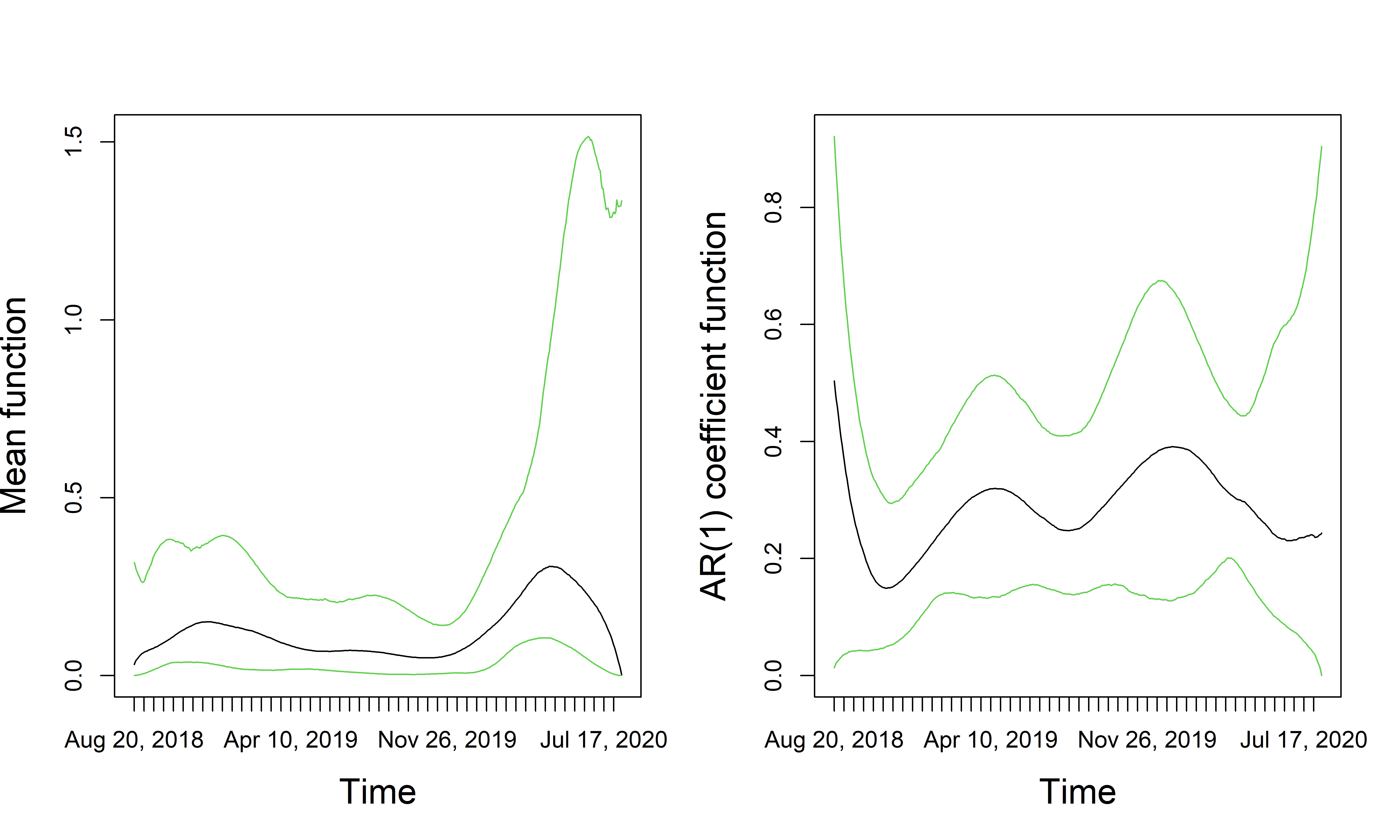}}
        \subfigure[$n=1000$]{\label{fig:c.63}\includegraphics[width=100mm, height=40 mm]{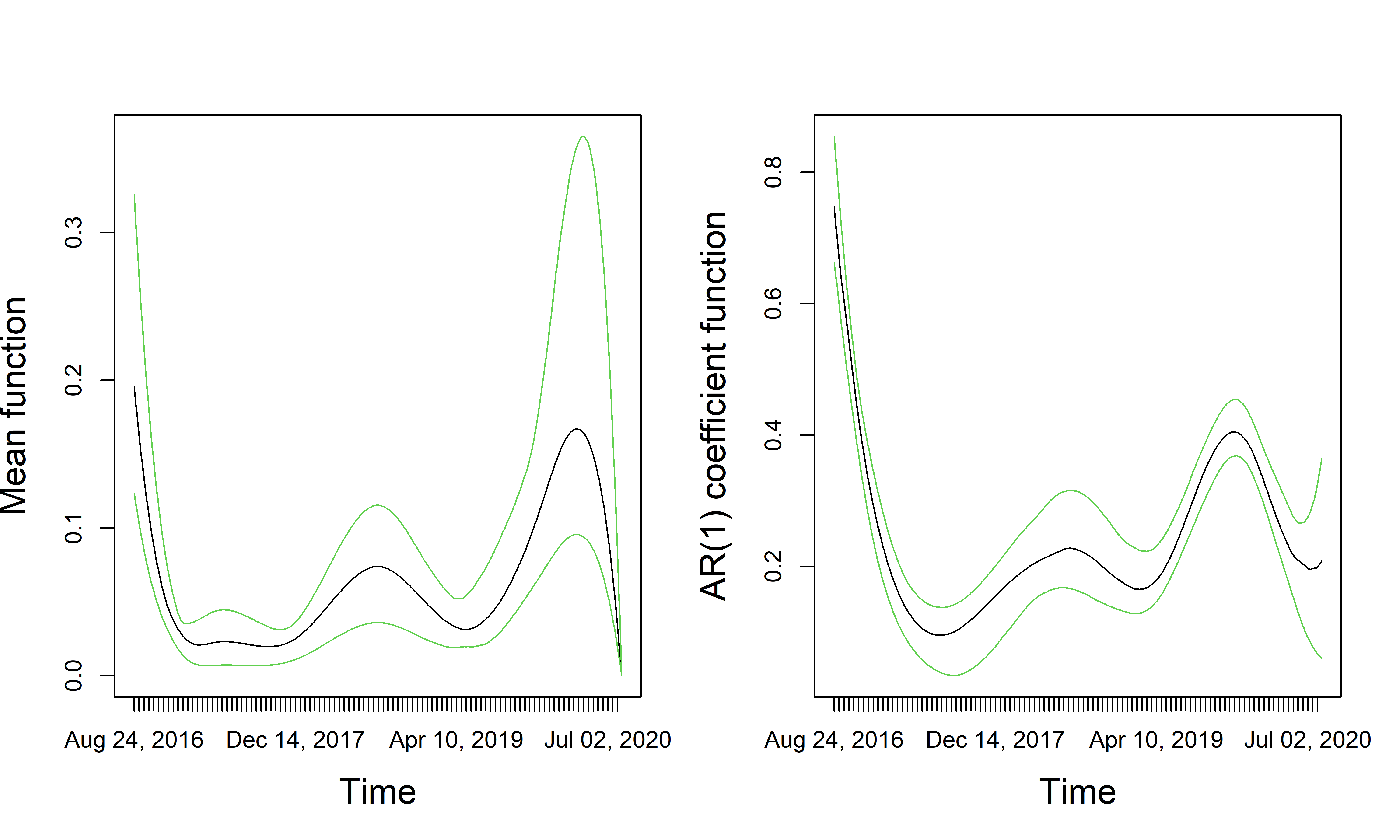}}
        \caption{NASDAQ data (tviGARCH(1,1) model) Estimated curve(black) along with the 95\% pointwise credible bands (green) are shown for T=200,500,1000 from top to bottom} 
        \label{fig:igarchdata}
\end{figure}

In Figure \ref{fig:garchdata} the sum of estimated coefficient functions $a(\cdot)+b(\cdot)$ is close to 1 for a significant time-horizon. This motivates us to also fit tviGARCH(1,1) to analyze the same NASDAQ data. The estimated functions are presented in  Figure \ref{fig:igarchdata} for the last $n=200, 500$ and $1000$ days. Table \ref{tab:AMSEtviGARCH1data} compares the AMSE scores for the same three methods as before with varying sample sizes. The estimated mean and AR(1) functions of Figure~\ref{fig:igarchdata} change a little from the estimated functions of tvGARCH(1,1) fit in Figure~\ref{AMSEtvGARCH1data}. Moreover, the effect of the three sample sizes is clear here with $n=1000$ showing very precise bands and can reveal an interesting time-varying pattern.

In terms of AMSE, one can see in Table~\ref{tab:AMSEtviGARCH1data} that the frequentist methods did worse than even the time-constant versions. The time-constant estimates were computed using the \texttt{rugarch} package in R. The Bayesian tviGARCH method provides significantly better AMSE uniform overall sample sizes. Here the mean function also shows higher volatility around the time when the pandemic struck us. Volatility due to the presidential election in 2016 can also be observed here.

\begin{table}[ht]
\centering
\caption{AMSE comparison: tviGARCH(1,1) model- Nasdaq data}
\begin{tabular}{rrrr}
  \hline
 & iGARCH(1,1) &  Frequentist tviGARCH(1,1) & Bayesian tviGARCH(1,1) \\ 
  \hline
    $n=200$ & 217.4988 & 278.4635 & \textbf{206.6886} \\  
$n=500$ & 96.5001 & 132.544 & \textbf{90.4456} \\ 
 $n=1000$ & 54.1171 & 260.4696 & \textbf{46.3704} \\ 
   \hline
\end{tabular}
\label{tab:AMSEtviGARCH1data}
\end{table}

\subsection{Model comparison}\label{ssc:modelcomp}
For the analysis of NASDAQ data, we have used two different models and thus it is pertinent to answer how should one choose between a competing class of models. We provide some measures in this subsection to decide between these two competing models. We start by comparing the performance of tvGARCH and tviGARCH models in terms of Bayes factor \citep{kass1995bayes}. Our calculation of the Bayes factor is based on the posterior samples using the harmonic mean identity of \cite{neton1994approximate}. Let us denote $B_{200}, B_{500}$ and $B_{1000}$ as the Bayes factors for three sample sizes where $$B_{i}=\frac{P(D^{(i)}|\textrm{tvGARCH})}{P(D^{(i)})|\textrm{tviGARCH})}$$ for sample size $i$ and the corresponding dataset $D^{(i)}$. The values we obtain are $2\log(B_{200})=8.16$, $2\log(B_{500})=19.08$ and $2\log(B_{1000})=24.14$. According to guidelines from section 3.2 of \cite{kass1995bayes}, there is `positive' evidence in favor of tvGARCH for sample sizes 200 and 500. However, the same evidence becomes `strong' for sample size 1000. 

We also try to address out-of-sample predictive performance comparison here. Note that for a time-varying GARCH or time-varying iGARCH model this is generally a difficult task due to the assumed non-stationarity of the model. Thus, we take following approach to calculate out of sample joint predictive log-likelihoods for model comparison. Let us assume we have the data $D^{(n)}$ with $n$ data points. To evaluate the joint predictive log-likelihood for the last $m (< n)$ most recent data points, we fit the models in~\eqref{modelgarch} and \eqref{modeligarch} with the first $n-m$ data. Note that the assumed time horizons for these two models are $n$. Based on the estimated B-spline coefficients, and other parameters from each model, we can compute the joint predictive log-likelihood of the last $m$ data points as $$L_m^{(n)}=\frac{1}{m}\sum_{i=n-m+1}^{n}\frac{1}{2}\left\{-X_{i}^2/\hat{\sigma}_{i}^2-\log(\hat{\sigma}_{i})-\log(2\pi)\right\},$$ where $\hat{\sigma}_{i}^2=\hat{\mu}(i/n)+\hat{a}_1(i/n)X_{i-1}^2+\hat{b}_1(i/n)\hat{\sigma}_{i-1}^2$. For the tviGARCH model, we have $\hat{b}_1(\cdot)=1-\hat{a}_{1}(\cdot)$.

Using this predictive log-likelihood we decide to evaluate the two fits from tvGARCH and tviGARCH in the following manner. For each of the sample sizes, we run it on three separate regimes of the data, the full data, and two halves of the data. In all these 8 settings, ( three sample sizes, three possible regions of the data, but the latest half of the 1000-sized data is the same as the full data for sample size 500) we compute 10, 20, and 50 steps ahead forecast. We tabulate these results in Table \ref{tab:modelcomparison}. One can see that generally speaking, there is somewhat conclusive evidence towards the iGARCH model for a smaller sample size. This supports our motivation why we additionally provide a tviGARCH(1,1) modeling on the same dataset.

Based on our model comparison exercises, we have an interesting phenomenon where for in-sample model fit, tvGARCH is better but in terms of out-of-sample prediction, tviGARCH outperforms tvGARCH in almost all the cases. Note that, tvGARCH has one additional free parameter and thus is expected to fit the data better but since the estimated $a_1(\cdot)$ and $b_1(\cdot)$ coefficients are close to one satisfying the iGARCH formulation, the out-of-sample performance for tviGARCH may have exceeded that for tvGARCH.

\begin{table}[!ht]
\caption{Joint log-likelihood for 10, 20, 50 steps ahead: Comparing tvGARCH(1,1)/tviGARCH(1,1) model- NASDAQ data. Better model is in bold}
\centering
\label{tab:modelcomparison}
\begin{tabular}{| l | l | | c | c || c | c || c | c ||}
\hline
& Steps & \multicolumn{2}{c||}{Full} &  \multicolumn{2}{c||}{First Half} & \multicolumn{2}{c||}{Second (Latest) Half}   \\
 $n$ & ($m$) & GARCH & iGARCH & GARCH & iGARCH & GARCH & iGARCH \\
\hline\hline
200 & 10 & $-2.1\times 10^{8}$ & \textbf{$-$1.982} & $-$3.532 & \textbf{$-$3.154} & $-8.9\times 10^{6}$ & \textbf{$-$2.251}\\ 
& 20 & $-$2.689 & \textbf{$-$1.881} & $-$14.889 & \textbf{$-$2.475} & $-$90281 & \textbf{$-$2.278}\\
& 50 & $-$3.640 & \textbf{$-$2.487} & $-$86.81 & \textbf{$-$2.877} & $-$5.839 & \textbf{$-$4.221}\\
\hline 

500 & 10 & $-$2.842 & \textbf{$-$2.068} & $-$10.260 & \textbf{$-$1.898} & $-$1161 & \textbf{$-$2.079}\\ 
 & 20 & $-$2.341 & \textbf{$-$1.848} & $-$3897 & \textbf{$-$1.407} & $-$343.49 & \textbf{$-$1.893} \\
 & 50 & $-$2.147 & \textbf{$-$2.112} & $-$44.381 & \textbf{$-$3.061} & $-$932.7 & \textbf{$-$2.371}\\ \hline

1000 & 10 & \textbf{$-$1.856} & $-$1.936 & \textbf{$-$2.499} & $-$2.790 & $-$2.842 & \textbf{$-$2.068}\\ 
 & 20 & $-$1.911 & \textbf{$-$1.789} & $-$1.999 & \textbf{$-$1.903} & $-$2.341 & \textbf{$-$1.848}\\
 & 50 & $-$2.266 & \textbf{$-$2.214} & $-$1.893 & \textbf{$-$1.536} & $-$2.147 & \textbf{$-$2.112}\\ \hline

\end{tabular}
\end{table}

As per the suggestion from a reviewer, we also add a one-step-ahead point forecasting exercise between these models. Here the computation method remains the same as outlined in the predictive log-likelihood computation however we only restrict ourselves to $m=1$ to make the discussion concise. For this part of the exercise we choose to compute posterior mean of $(X_n^2-\hat{\sigma}^2)^2$ where to ensure out-of-sample prediction $\hat{\sigma}^2$ is estimated solely based on $X_1,\ldots,X_{n-1}$. As one-step-ahead forecasts can be prohibitively misleading given it depends so much on one single location, we decide to take an average of over 15 random starting points over the entire time spectrum of 10 years resulting in 2518 points). For each of the sample sizes, we tabulate the performance in the following Table \ref{tab:garchigarchpoint}. Note that, here we are only comparing the two Bayesian time-varying models to see which one fits our data better. The advantage of predicting the future coefficients using B-spline is not available in the kernel-based frequentist method and thus is not included here in the discussion.

\begin{table}[ht]
\centering
\caption{One-step-ahead out-of-sample forecast for NASDAQ data: Comparing tvGARCH(1,1)/tviGARCH(1,1) model}
\begin{tabular}{rrr}
  \hline
 &  Bayesian tvGARCH(1,1) &  Bayesian tviGARCH(1,1) \\ 
  \hline
    $n=200$ & 0.9914 & \textbf{0.3477}  \\  
$n=500$ & 1.5208  & \textbf{1.4102}  \\ 
 $n=1000$ & \textbf{1.6378} & 1.7033 \\ 
   \hline
\end{tabular}
\label{tab:garchigarchpoint}
\end{table}

\noindent One can see, we again observe the same advantage of tviGARCH modeling over the tvGARCH one for smaller sample sizes. This is an interesting find of this paper in the context of the Bayesian model fitting to these datasets.

\section{Discussion and Conclusion}
	\label{sec:discussion}
In this paper, we consider a Bayesian estimation framework for time-varying conditional heteroscedastic models. Our prior specifications are amenable to Hamiltonian Monte Carlo for efficient computation. One of the key motivations towards going to a Bayesian regime was to achieve reasonable estimation even for a small sample size. Our simulation coverage shows good performance for all three models tvARCH, tvGARCH, tviGARCH for both small and large sample sizes. Importantly, in all three of the cases, we were able to establish posterior contraction rates. These calculations are, to the best of our knowledge, the first such work in even simple dependent models let alone the complicated recursions that these conditional heteroscedastic models demand. Moreover, the assumptions on the true functions and the number of moments needed were very minimal. An interesting future theoretical work would be to calculate posterior contraction rate with respect to empirical $\ell_2$-distance which is a more desirable metric for function estimation. While analyzing real data, we see that the parameter curves vary significantly for the intercept terms, but not that much for AR or CH parameters. The associated codes to fit the three models are available at \url{https://github.com/royarkaprava/tvARCH-tvGARCH-tviGARCH}.

As future work, it will be interesting to explore multivariate time-varying volatility modeling \citep{multigarch,multigarch2} through a Bayesian framework similar to ours. Another interesting time-heterogeneity that we plan to explore through the glass of a Bayesian framework is regime-switching CH models where instead of the smooth time-varying functions the parameters change abruptly. We have a brief discussion in section \ref{ssc:modelcomp} about how to choose between competing models. Those discussions can easily be extended to choose a proper number of lags or to choose between different types of ARCH/GARCH models. We believe this would provide an interesting parallel to the usual penalized likelihood-based methods for model selection in time-series. Finally note that, even though we do provide some insights onto future prediction for these datasets for real data applications, that was not the main focus in this paper. Forecasting for the time-varying model is extremely tricky since it requires `estimation' of the future parameter values. Although in-filled asymptotics can help in this regard, still the literature so far is very sparse in this direction for both Bayesian and frequentist regimes. We plan to explore this extensively in near future.

\section*{Acknowledgement}
We would like to thank the editor, the associate editor, and anonymous referees for their constructive suggestions that improved the quality of the manuscript.

	\bibliographystyle{apalike}
	\bibliography{main}

\newpage

\section{Proof of Theorems}\label{sec:proof}
 We study the frequentist property of the posterior distribution in increasing $n$ regime assuming that the observations are coming from a true density $f_0$ characterized by the parameter $\kappa_0$. We follow the general theory of \cite{ghosal2000convergence} to study posterior contraction rate for our problem. In Bayesian framework, the density $f$ is itself a random measure and has distribution $\Pi$ which is the prior distribution induced by the assumed prior distribution on $\kappa$. The posterior distribution of a neighborhood $U_n=\{f:d(f,f_0)<\epsilon_n\}$ around $f_0$ given the observation $X^{(n)}=\{X_0,X_1,\ldots,X_n\}$ is
$$
\Pi_n(U_n^c|X^{(n)})=\frac{\int_{U_n^c}f(X^{(n)})d \Pi(\kappa)}{\int f(X^{(n)})d\Pi(\kappa)}
$$

\subsection{General proof strategy}\label{sec:genproof}

The posterior consistency would hold if above posterior probability almost surely goes to zero in $F_{\kappa_0}^{(n)}$ probability as $n$ goes to $\infty$, where $F_{\kappa_0}^{(n)}$ is the true distribution of $X^{(n)}$. Recall the definition of posterior contraction rate; for a sequence $\epsilon_n$ if $\Pi_n(d(f,f_0)|X^{(n)}\geq M_n\epsilon_n|X^{(n)})\rightarrow 0$ in $F_{\kappa_0}^{(n)}$-probability for every sequence $M_n\rightarrow\infty$, then the sequence $\epsilon_n$ is called the posterior contraction rate. If the assertion is true for a constant $M_n=M$, then the corresponding contraction rate becomes slightly stronger. 

Note that for two densities $f_{0}, f$ characterized by $\kappa_0$ and $\kappa$ respectively, the Kullback-Leibler divergences are given by 
\begin{gather*}
KL(\kappa_0, \kappa) = \int f_0\log{\frac{f_0}{f}}=E_{\kappa_0}\left[\log\frac{\IP_{Q_{\kappa_0}}(X_0)\prod_{i=1}^n\IP_{\kappa_0}(X_i|\sF_{i-1},\lambda_0)}{\IP_{Q_{\kappa}}(X_0)\prod_{i=1}^n\IP_{\kappa}(X_i|\sF_{i-1},\lambda_0)}\right].
\end{gather*}

\noindent Assume that there exists a sieve in parameter space such that $\Pi(W_n^c)\leq \exp(-(C_n+2)n\epsilon_n^2)$ and we have tests $\chi_n$ such that $$\IE_{\kappa_0}(\chi_n)\leq e^{-L_nn\epsilon_n^2/2}\quad \sup_{\kappa\in W_n: d^2(f,f_0)>L_{n}\epsilon_n^2}\IE_{\kappa}(1-\chi_n)\lesssim e^{-L_{n}n\epsilon_n^2}$$ for some $L_{n}>C_n+2$. Say $U_n=\{f:d^2(f,f_0)\leq L_n\epsilon_n^2\}$  and $S_n=\{\int \frac{f(X^n)}{f_0(X^n)}d\Pi(\kappa)\geq\Pi_{n}(\frac{1}{n}KL(\kappa_0,\kappa)<\epsilon_n)\exp(-C_{n}n\epsilon_n^2)\}$. We can bound the posterior probability from above by,
\begin{align}
    \Pi_n(d(f,f_0)\geq M_n\epsilon_n|X^{(n)}) &\leq \chi_n+(1-\chi_n)\frac{\int_{U_n^c}f(X^n)d \Pi(\kappa)}{\int f(X^{(n)})d\Pi(\kappa)}\nonumber\\
    &=\chi_n+(1-\chi_n)\frac{\int_{U_n^c}\frac{f(X^{(n)})}{f_0(X^{(n)})}d \Pi(\kappa)}{\int \frac{f(X^{(n)})}{f_0(X^{(n)})}d\Pi(\kappa)}\nonumber\\
    &\leq \chi_n+\mathbb{1}\{S_n^c\} + (1-\chi_n)\frac{\int_{U_n^c}\frac{f(X^{(n)})}{f_0(X^{(n)})}d \Pi(\kappa)}{\exp(-C_nn\epsilon_n^2)\Pi_{n}\{\frac{1}{n}KL(\kappa_0,\kappa)<\epsilon_n\}} \nonumber\\
    &\leq \chi_n+\mathbb{1}\{S_n^c\} + \frac{\exp(C_nn\epsilon_n^2)}{\Pi_{n}\{\frac{1}{n}KL(\kappa_0,\kappa)<\epsilon_n\}}(1-\chi_n)\frac{\int_{U_n^c}f(X^{(n)}}{f_0(X^{(n)})}d \Pi(\kappa)
\end{align}
Taking expectation with respect to $\kappa_0$, first term goes to zero by construction of $\chi_n$. The second term $\IE_{\kappa_0}\mathbb{1}\{S_n^c\}$ goes to zero due to Lemma 8.21 of \cite{ghosal2017fundamentals} for any sequence $C_n\rightarrow \infty$. We would require that $\Pi_{n}\{\frac{1}{n}KL(\kappa_0,\kappa)< \epsilon_n\} \geq \exp(-n\epsilon_n^2)$. Then for the third term,
\begin{align}
    &\IE_{\kappa_0}\exp((C_n+1)n\epsilon_n^2)(1-\chi_n)\frac{\int_{U_n^c}f(X^{(n)})}{f_0(X^{(n)})}d \Pi(\kappa)=\exp((C_n+1)n\epsilon_n^2)\int_{U_n^c}f(X^{(n)})(1-\chi_n)d\Pi(\kappa)\nonumber\\
    &\quad\leq\exp(C_n+1)n\epsilon_n^2)\left[\int_{U_n^c\cap W_n}f(X^{(n)})(1-\chi_n)d\Pi(\kappa) + \Pi(W_n^c)\right]\nonumber\\&\quad=\exp((C_n+1)n\epsilon_n^2)\left[\sup_{\kappa\in W_n: d^2(f,f_0)>L_{n}\epsilon_n^2}\IE_{\kappa}(1-\chi_n) + \Pi(W_n^c)\right]\lesssim \exp(-n\epsilon_n^2).
\end{align}

\noindent Thus we adhere to the following plan
\begin{plan}\label{req}
The proof had three major parts as follows:
\begin{itemize}
    \item[(i)](Prior mass Condition) We need $\Pi_{n}\{\frac{1}{n}KL(\kappa_0,\kappa)<\epsilon_n\} \geq \exp(-n\epsilon_n^2)$, 
    \item[(ii)](Sieve) Construct the sieve $W_n$ such that $\Pi(W_n^c)\leq \exp(-(C_n+2)n\epsilon_n^2)$ and
    \item[(iii)](Test construction) Construct exponentially consistent tests $\chi_n$.
\end{itemize}
\end{plan}

We first study the contraction properties with respect to $d^2(f,f_0)=r_n^2(f, f_0)=-\frac{1}{n}\log\int\sqrt{ff_0}$ and then show that the same rate holds for average Hellinger $\frac{1}{n}d_H^2(f,f_0)$. Note that $L_n$ can be taken as $L_n=M_n^2$.

\subsection{Proof of Theorem 1}\label{sec:proofthm1}
For the sake of technical convenience we show our proof by fixing lag order $p$ at $p=1$, however the results are easily generalizable for any fixed order $p$. All the proofs go through for higher lags with the same technical tools.
\subsubsection{KL Support}\label{sec:klbound}
The likelihood based on the parameter space $\kappa$ is given by
$
\IP_{\kappa}(X_0)\prod_{i=1}^n\IP_{\kappa}(X_i|X_{i-1}).
$ Let $Q_{\kappa, t}(X_i)$ be the distribution of $X_i$ with parameter space $\kappa$.

\noindent We have 
\begin{eqnarray}\label{eq:def R}
    R&=&\frac{1}{2}\log\frac{\prod_{i=1}^n\IP_{\kappa_0}(X_i|\sF_{i-1},\sigma_0)}{\prod_{i=1}^n\IP_{\kappa}(X_i|\sF_{i-1},\sigma_0)}\nonumber\nonumber
    \\&=&\frac{1}{2}\sum_{i=1}^n \left[-\log (\mu_0(i/n)+a_{01}(i/n)X_{i-1}^2)+\log (\mu(i/n)+a_{1}(i/n)X_{i-1}^2)\right] \nonumber \\ \nonumber 
    && \qquad \qquad +\sum_{i=1}^n\left[ X_i^2\left(\frac{1}{\mu_0(i/n)+a_{01}(i/n)X_{i-1}^2}-\frac{1}{\mu(i/n)+a_{1}(i/n)X_{i-1}^2}\right) \right] \\ 
    &=& \frac{1}{2}\sum_{i=1}^n\{(\log \sigma_i^2-\log \sigma_{0i}^2)+X_i^2(1/\sigma_{0i}^2- 1/\sigma_i^2)\}= \frac{1}{2}(I+II) \text{ (say)} \label{eq:IandII}
\end{eqnarray}
where $\sigma_i^2=\mu(i/n)+a_1(i/n)X_{i-1}^2$ and $\sigma_{0i}^2=\mu_0(i/n)+a_{01}(i/n)X_{i-1}^2$. Then
\begin{eqnarray}\label{eq:KLandV}
KL(\kappa^n_0, \kappa^n)=\IE_{\kappa_0}(R).
\end{eqnarray} Our first goal is to estimate these two quantities in terms of the distances between $\mu(\cdot), \mu_0(\cdot)$ and $a_1(\cdot), a_{01}(\cdot)$. In the light of $\IE_{\kappa_0}(I+II)=\IE_{\kappa_0}(I)+\IE_{\kappa_0}(II)$,  we bound the expectation of $I$ and $II$ separately. 

\medskip
\noindent \underline{Bounding the first term I:} Note that, by the means of mean value theorem, for a random variable 
$\sigma_{*t}^2$ between $\sigma_i^2$ and $\sigma_{0i}^2$,
\begin{eqnarray}
\log \sigma_i^2- \log \sigma_{0i}^2&=& \frac{1}{\sigma_{*t}^2}(\mu(i/n)-\mu_0(i/n)+(a_{1}(i/n)-a_{01}(i/n)) X_{i-1}^2) \nonumber \\
|\log \sigma_i^2- \log \sigma_{0i}^2|&\leq& \frac{1}{\rho}(|\mu(i/n)-\mu_0(i/n)|+|a_{1}(i/n)-a_{01}(i/n)|) \label{eq:Rfirstterm}
\end{eqnarray}
where $\sigma_{*i}^2>\rho$ for the first term and $\sigma_{*i}^2>\rho X_{i-1}^2$ is used for the second term due to the assumption (A.2). Thus the first term $I$ satisfies

$$\IE_{\kappa_0}(I) \leq \IE_{\kappa_0}(|I|) \lesssim T(\|\mu-\mu_0\|_{\infty} +\|a_1-a_{01}\|_{\infty}),$$ in the light of assumption (A.3). 

\noindent \underline{Bounding the second term II:}
For the second term we proceed as follows:
$$X_i^2(1/\sigma_{0i}^2- 1/\sigma_i^2)\} \leq \frac{1}{\sigma_i^2\sigma_{0i}^2}|X_i^2\sigma_{0i}^2-X_i^2\sigma_i^2|\leq \frac{X_i^2}{\sigma_{0i}^2}\left(\frac{|\mu(i/n)-\mu_0(i/n)|}{\mu(i/n)}+\frac{|a_1(i/n)- a_{01}(i/n)|X_{i-1}^2}{a(i/n)X_{i-1}^2}\right)$$

\noindent where we use the fact that $\mu(i/n)>\rho$ and $a(i/n)>\rho$ due to closeness of $(\mu,a)$ and $(\mu_0,a_0)$. Consequently, we have the deterministic inequality

$$X_i^2(1/\sigma_{0i}^2- 1/\sigma_i^2)\}\leq \frac{X_i^2}{\sigma_{0i}^2}\frac{\|\mu-\mu_0\|_{\infty}+\|a_1-a_{01}\|_{\infty}}{\rho}.$$

\noindent Taking expectation under truth we arrive at similar upper bound for  $\IE_{\kappa_0}(II)$ as $\IE_{\kappa_0}(I)$ since $$\IE(X_i^2/\sigma_{0i}^2)=\IE(\IE(X_i^2/\sigma_{0i}^2|X_{i-1}))=1.$$

\ignore{
\begin{eqnarray}\label{eq:Rsecondterm}
X_i^2(1/\sigma_{0i}^2- 1/\sigma_i^2)\} &\leq& \frac{1}{\sigma_i^2\sigma_{0i}^2}|X_i^2\sigma_{0i}^2-X_i^2\sigma_i^2| \nonumber\\
&\leq& \frac{1}{\rho}\left(X_i^2\frac{|\mu(i/n)-\mu_0(i/n)|}{\mu_{0}(i/n)}+X_i^2\frac{|a_1(i/n)- a_{01}(i/n)|X_{i-1}^2}{a_{01}(i/n)X_{i-1}^2}\right)\nonumber \\
&\leq& \frac{X_i^2}{\rho^2}(\|\mu-\mu_0\|_{\infty}+\|a_1-a_{01}\|_{\infty})
\end{eqnarray} by assumption (A.1-A.2). Thus the second term $II$ satisfies

$$\IE_{\kappa_0}(II) \leq \IE_{\kappa_0}(|II|) \lesssim T(\|\mu-\mu_0\|_{\infty} +\|a_1-a_{01}\|_{\infty}),$$ in the light of finite second moment from assumption (A.3). 
}

\noindent Thus we conclude that,
\begin{align}
   \frac{1}{n}\IE(R)\lesssim \|\mu-\mu_0\|_{\infty}+\|a_1-a_{01}\|_{\infty}.\label{eq:KLineq1}
\end{align}

\subsubsection{Posterior contraction in terms of average negative log-affinity}
In this section, we focus on the requirements to calculate posterior contraction rate as in Section~\ref{req}.
 We first show posterior consistency in terms of average negative log-affinity which is defined as $r_n^2(f_1,f_2)=-\frac{1}{n}\log\int f_1^{1/2}f_2^{1/2}$ between $f_1$ and $f_2$. Here, we have $f_1=\prod_{i=1}^n P_{\kappa_1}(X_{i}|X_{i-1})$. Then we show that, having $r_n^2(f_1,f_0)\lesssim\epsilon_n^2$ implies that our distance metric $d_{2,n}^2(f_1,f_0)\lesssim \epsilon_n^2$.

\noindent Proceeding with the rest of the proof of Theorem 1, we use the results of B-Splines,
$\|\mu-\mu_0\|_{\infty}\leq \|\alpha-\alpha_0\|_{\infty},$ where $\alpha = \{\alpha_j\}$ and $\|a_1-a_{01}\|_{\infty}\leq \|\gamma-\gamma_{0}\|_{\infty}$, where $\gamma_{j}=\theta_{1j}M_1,$ such that $\gamma_{j}<1$. The H\"older smooth functions with regularity $\iota$ can be approximated uniformly up to order $K^{-\iota}$ with $K$ many B-splines. Thus we have $\epsilon_n\gtrsim\max\{K_{1n}^{-\iota_1},K_{2n}^{-\iota_2}\}$. 

We start by providing a lower bound of the prior probability as required by (i) in Plan \ref{req}. We also have the result~\eqref{eq:KLineq1} and the prior probabilities
$\Pi(\|\alpha-\alpha_0\|_{\infty} \lesssim \epsilon_n, \|\gamma-\gamma_{0}\|_{\infty} \lesssim \epsilon_n) \gtrsim \epsilon_n^{K_{1n}+K_{2n}}$ based on the discussion of A2 from \cite{shen2015adaptive}. The rate of contraction cannot be better than the parametric rate $n^{-1/2}$, and so $\log (1/\epsilon_n)\lesssim \log T$. Thus (i) requires that in terms of pre-rate $\bar{\epsilon}_n$ as $(K_{1n}+K_{2n})\log n\lesssim n\bar{\epsilon}_n^2$. Now we construct a sequence of test functions having exponentially decaying Type I and Type II error. In our problem, we consider following sieve 

\begin{align}\label{eq:sieve}
    W_n&=\{K_1,K_2,\alpha,\gamma:K_1\leq K_{1n},K_2\leq K_{2n}, \|\alpha\|_{\infty}\leq A_{n}, \min(\alpha,\gamma)>\rho_n,  \nonumber\\&\qquad\qquad\qquad\qquad\qquad\qquad\qquad\qquad\gamma\leq 1-A_n/B_n,\sigma_0^2\leq B_n, A_n<B_n \},
\end{align}
where $A_n, B_n$ are at least polynomial in $n$ and $\sigma_0^2$ is the standard deviation of $X_0^2$ and $K_n=\max\{K_{1n},K_{2n}\}$. We take $\rho_n=O(n^{-a})$ with $a<1$, $A_n=O(n^{a_1}), B_n=O(n^{a_2})$ with $a_2>a_1$ for technical need. 

 We need to choose these bounds carefully so that we have $\Pi(W_n^c)\leq \exp(-(1+C_1)n\epsilon_n^2)$, which depend on tail properties of the prior. We have, 
\begin{eqnarray*}
\Pi(W_n^c)&\leq& \Pi(K_1>K_{1n})+\Pi(K_2>K_{2n})+\Pi\{\alpha_{K_{1n}}\notin[\rho_n, A_n]^{K_{1n}}\}+\Pi\{\gamma_{K_{2n}}\notin\left[\rho_n, 1-\frac{A_n}{B_n}\right]^{K_{2n}}\}\\&\quad&+\Pi\{\lambda^0>B_n\}
\end{eqnarray*}
where $\alpha_{K_{1n}}$ is the vector of full set of coefficients of length $K_{1n}$ and $\gamma_{K_{2n}}$ is the vector of coefficients of length $K_{2n}$. The quantity $\Pi[\alpha_{K_{1n}}\notin[\rho_n), A_n]^{K_{1n}}$ can be further upper bounded by $K_{1n}\Pi(\alpha_{1}\notin[\rho_n, A_n)])\leq K_{1n}\exp\{-R_1n^{a_3}\}$, for some constant $R_1, a_3>0$ which can be verified from the discussion of the assumption A.2 of \cite{shen2015adaptive} for our choice of exponential prior. On the other hand, $$\Pi\{\gamma_{K_{2n}}\notin\left[\rho_n, 1-\frac{A_n}{B_n}\right]^{K_{2n}}\}\leq K_{2n}\Pi(\gamma_1\notin[\rho_n, 1-\frac{A_n}{B_n}])\leq K_{2n}\exp\{-R_2n^{a_4}\},$$ for some constant $R_2, a_4>0$ which can be verified from the proof of \cite{roy2018high}. 
Hence, $\Pi(W_n^c)\lesssim F_1(K_{1n})+F_2(K_{2n})+(K_{1n}+K_{2n})\exp\{-Rn^{a_5}\}$.  The two functions $F_1$ and $F_2$ in the last expression stand for the tail probabilities of the prior of $K_1$ and $K_2$. We can calculate their asymptotic order as, $F_1(x)=\Pi(K_{1}>x)\asymp\exp\{-x(\log x)^{b_{13}}\}$ and $F_2(x)=\Pi(K_{2}>x)\asymp\exp\{-x(\log x)^{b_{23}}\}$. We need $\Pi(W_n^c)\lesssim\exp\{-(1+C_n)n\epsilon_n^2\}$. Hence, we calculate pre-rate from the following equation for some sequence $H_n\rightarrow\infty$,

\begin{align}
\label{eq:rate sieve}
{K}_{1n} (\log n)^{b_{13}}+{K}_{2n}(\log n)^{b_{23}}\gtrsim H_nn\bar\epsilon_n^2, \quad \log(K_{1n}+K_{2n})+H_nn\bar\epsilon_n^2\lesssim n^{a_5}.
\end{align}

\noindent Now, we construct test $\chi_n$ such that $$\IE_{\kappa_0}(\chi_n)\leq e^{-L_nn\epsilon_n^2/2}\quad \sup_{\kappa\in W_n: r_n^2(\kappa,\kappa_0)>L_n\epsilon_n^2}\IE_{\kappa}(1-\chi_n)\lesssim e^{-L_nn\epsilon_n^2}$$ for some $L_n>C_n+2$.  To construct the test, we first construct the for point alternative $H_0:\kappa=\kappa_0$ vs $H_1:\kappa=\kappa_1$. The most powerful test for such problem is Neyman-Pearson test $\phi_{1n}=\mathbb{1}\{f_1/f_0\geq 1\}$.  For $r_n^2> L_n\epsilon_n^2$, we have

$$\IE_{\kappa_0}\phi_{1n}=\IE_{\kappa_0}(\sqrt{f_1/f_0}\geq 1)\leq \int \sqrt{f_1f_0}\leq \exp(-L_nn\epsilon_n^2),$$ 

$$\IE_{\kappa_1}(1-\phi_{1n})=\IE_{\kappa_1}(\sqrt{f_0/f_1}\geq 1)\leq \int \sqrt{f_0f_1}\leq \exp(-L_nn\epsilon_n^2).$$

\noindent It is natural to have a neighborhood around $\kappa_1$ such the Type II error remains exponentially small for all the alternatives in that neighborhood under the test function $\phi_{1n}$. By Cauchy-Schwarz inequality, we can write that $$\IE_{\kappa}(1-\phi_{1n})\leq \{\IE_{\kappa_1}(1-\phi_{1n})\}^{1/2}\{\IE_{\kappa_1}(f/f_1)^2\}^{1/2}.$$In the above expression, the first factor is already exponentially decaying. The second factor can be allowed to grow at most of order $e^{cn\epsilon_n^2}$ for some positive small constant $c$. We, in fact establish that the second factor is bounded by our choice of the sieve through the following claim.

\begin{claim}\label{claim1}
 $\IE_{\kappa_1}(f/f_1)^2$ is bounded for every $\kappa$ such that
\begin{eqnarray}\label{eq:condition1}
r_1, r_2 \leq \frac{\rho_n}{4n}, \text{ where }\|\mu-\mu_1\|_{\infty}=r_1, \|a-a_1\|_{\infty}=r_2.
\end{eqnarray}
\end{claim}

\begin{proof}
\noindent We have, in the light of AM-GM inequality, 

\begin{eqnarray}\label{eq:crucial}
\IE_{\kappa_1}(f/f_1)^2=\int\frac{f^2}{f_1^2}f_1=\IE_{f}\frac{f}{f_1}=\IE_{f}\prod_{i=1}^n\frac{f(X_i|X_{i-1})}{f_1(X_i|X_{i-1})}\leq\frac{1}{n}\sum_{i=1}^n\IE_{f}\left[\left(\frac{f(X_i|X_{i-1})}{f_1(X_i|X_{i-1})}\right)^n\right]
\end{eqnarray}

\noindent Towards uniformly bounding the summand, we first compute and provide a deterministic upper bound to the conditional expectation of the same given $X_{i-1}$. Denote $\sigma^2=\mu(i/n)+a_1(i/n)X_{i-1}^2$ and $\sigma_1^2=\mu_1(i/n)+a_{11}(i/n)X_{i-1}^2$. Note that if $(n+1/2)\sigma_1^2-n\sigma^2 \geq 0$, then
\begin{align}
     \nonumber \IE_{f}\left[\left(\frac{f(X_i|X_{i-1})}{f_1(X_i|X_{i-1})} \right)^n | X_{i-1} \right] &= \int_{x}\frac{\{f(X_i=x|X_{i-1})\}^{n}}{\{f_1(X_i=x|X_{i-1})\}^n}f(X_i=x|X_{i-1})dx \\ &=\left( \frac{\sigma_1^2}{\sigma^2} \right)^{n/2}\frac{1}{\sqrt{2 \pi}\sigma}\int_{x}exp\left(-\frac{Tx^2}{2}(\frac{1}{\sigma^2}-\frac{1}{\sigma_1^2}) -\frac{x^2}{2\sigma^2}\right) dx \nonumber \\
    &=\left( \frac{\sigma_1^2}{\sigma^2} \right)^{n/2}\frac{\sigma_0}{\sigma} \text{ (where } \frac{1}{\sigma_0^2}=\frac{T+1}{\sigma^2}-\frac{n}{\sigma_1^2} )\nonumber \\
    &=\left( \frac{\sigma_1^2}{\sigma^2} \right)^{n/2}\frac{\sigma_1}{\sqrt{(n+1)\sigma_1^2-n\sigma^2}} \nonumber \\
    & \leq \sqrt{2} \left( \frac{\sigma_1^2}{\sigma^2} \right)^{n/2}\label{eq:tobdexp}.
\end{align}
Now since  
$$\frac{\sigma^2-\sigma_1^2}{\sigma_1^2} \leq \frac{|\sigma^2-\sigma_1^2|}{\sigma_1^2}\leq  \frac{|\mu_(i/n)-\mu_{1}(i/n)|}{\mu_1(i/n)}+\frac{|a_1(i/n)-a_{11}(i/n)|}{a_{11}(i/n)}\leq \frac{r_1+r_2}{\rho_n}$$ due to the definition of sieve at (\ref{eq:sieve}), we have in the light of assumption (\ref{eq:condition1}), $(n+1/2)\sigma_1^2-n\sigma^2 \geq 0$. Assumption (\ref{eq:condition1}) allows us to bound the summand in (\ref{eq:crucial}) before expectation as well. Note that, in the light of (\ref{eq:condition1}), for large $n$, we have the deterministic inequality
$$\left(\frac{\sigma_1^2}{\sigma^2} \right)^{n/2}=\left(1+\frac{\sigma_1^2-\sigma^2}{\sigma^2}\right)^{n/2} \leq \left(1+\frac{r_1+r_2}{\rho_n}\right)^{\frac{n}{2}}\leq \left(1+\frac{1}{2n}\right)^{\frac{n}{2}}\approx e^{1/4}.$$
\end{proof}

\noindent The test function $\chi_n$ satisfying exponentially decaying Type I and Type II probabilities is then obtained by taking maximum over all tests $\phi_{jn}$'s for each ball, having above radius. Take $\chi_n=\max_j\phi_{jn}$. Type I and Type II probabilities are given by $\IP_0(\chi_n)\leq \sum_{j}\IP_0\phi_{jn}\leq D_{n}\IP_0\phi_{jn}$ and $\sup_{\kappa\in W_n: r_n(\kappa,\kappa_0)>Mn\epsilon_n^2}\IP(1-\chi_n)\leq \exp(-cn\epsilon_n^2)$. Hence, we need to show that $\log D_{n}\lesssim n\epsilon_n^2$, where $D_n$ is the required number of balls of above radius needed to cover our sieve $W_n$. We have 
\begin{align}
   \log D_n&\leq \log D(r_1, \|\alpha\|_{\infty}\leq A_n, \min(\alpha)>\rho_n, \|\|_{\infty})+\log D(r_2, \|\gamma\|_{\infty}\leq 1-\frac{A_n}{B_n}, \min(\gamma)>\rho_n, \|\|_{\infty}) \nonumber\\
   &\leq K_{1n}\log(3K_{1n}A_n/r_1)+K_{2n}\log(3K_{2n}/r_2)
\end{align}
where $r_1,r_2$ are defined in (\ref{eq:condition1}). Given our choices of $A_n, B_n$ and $\rho_n$, the two radii $r_1$ and $r_2$ are some fractional polynomials in $n$. Thus $\log D_n\lesssim (K_{1n}+K_{2n})\log T$ which is required to be $\lesssim n\epsilon_n^2$ as in the prior mass condition. Based on~\eqref{eq:rate sieve}, we have 
$$K_{1n}\asymp n^{1/(2\iota_1+1)}(\log n)^{\iota/(2\iota+1)+(1-b_{13})/2}, K_{2n}\asymp n^{1/(2\iota_2+1)}(\log n)^{\iota/(2\iota_2+1)+1-b_{23}}$$ and a pre-rate $$\bar{\epsilon}_n=\max\bigg\{n^{-\iota_1/(2\iota_1+1)} (\log n)^{\iota_1/(2\iota_1+1)},n^{-\iota_2/(2\iota_2+1)} (\log n)^{\iota_2/(2\iota_2+1)}\bigg\}.$$ The actual rate will be slower that pre-rate. Now, the covering number condition, prior mass conditions and basis approximation give us $(K_{1n}+K_{2n})\log n\lesssim n\epsilon_n^2$ and $\epsilon_n\gtrsim\max\{K_{1n}^{-\iota_1},K_{2n}^{-\iota_2}\}$.
 Combining all these conditions, we calculate the posterior contraction rate $\epsilon_n$ equal to \begin{align*}
    \max\bigg\{&n^{-\iota_1/(2\iota_1+1)} (\log n)^{\iota_1/(2\iota_1+1)+(1-b_{13})/2},n^{-\iota_2/(2\iota_2+1)} (\log n)^{\iota_2/(2\iota_2+1)+(1-b_{23})/2}\bigg\}.
\end{align*}
\subsubsection{Posterior contraction in terms of average Hellinger}\label{ssc:equiv}
Write Reyni divergence as $$r_n^2=-\frac{1}{n}\log\int\sqrt{f_0f_1}=-\frac{1}{n}\log \IE_{\kappa_0}\sqrt{\frac{f_1}{f_0}}.$$We need to show $r_n^2\lesssim\epsilon_n^2$ implies that $d_{2,n}^2(\kappa_0,\kappa)\lesssim \epsilon_n^2$ as $\epsilon_n$ goes to zero. If $r_n^2\leq\epsilon_n^2$, we have $\left(\IE_{\kappa_0}\sqrt{\frac{f_1}{f_0}}\right)^{-1/n}\leq\exp(\epsilon_n^2)$ which implies for small $\epsilon_n^2$, we have $\left(\IE_{\kappa_0}\sqrt{\frac{f_1}{f_0}}\right)^{1/n}\geq 1-\epsilon_n^2$. By Cauchy-Squarz inequality $\left(\int\sqrt{f_0f_1}\right)^2\leq \int f_0\int f = 1$. Thus we have,

$$
1-\epsilon_n^2\leq \left(\IE_{\kappa_0}\sqrt{\frac{f_1}{f_0}}\right)^{1/n}\leq 1,
$$

\noindent Since $d_{H}^2(f_1,f_0)=2(1-\IE_{\kappa_0}\sqrt{\frac{f_1}{f_0}})$
$$\left(\IE_{\kappa_0}\sqrt{\frac{f_1}{f_0}}\right)^{1/n}=\left\{1-\left(1-\IE_{\kappa_0}\sqrt{\frac{f_1}{f_0}}\right)\right\}^{1/n}\approx1-\frac{1}{2n}d_{H}^2(f_1,f_0).$$ Thus $\frac{1}{n}d_{H}^2(f_1,f_0)\lesssim\epsilon_n^2$. Thus it is consistent under average Hellinger distance.

\subsection{Proof of Theorem 2}\label{sec:proofthm2}
Note that, the proof of Theorem \ref{thm:garch} follows via exactly same route. We just jot down the important differences from the proof of Theorem \ref{thm:arch}. For Theorem \ref{thm:garch} also, we restrict ourselves to tvGARCH(1,1) situation with the assurance that the proof easily extends to a general tvGARCH(p,q) case.

\subsubsection{KL support}
 Note that the KL support step from subsubsection \ref{sec:klbound} follows almost similarly with the modified 
$$\sigma_i^2=\mu(i/n)+a_1(i/n)X_{i-1}^2+b_1(i/n)\sigma_{i-1}^2 \text{  and }\sigma_{0i}^2=\mu_0(i/n)+a_{01}(i/n)X_{i-1}^2+b_{01}(i/n)\sigma_{0,i-1}^2$$
To our advantage, we also have the lower bound for the additional $b_1(\cdot)$ function from assumption B.2. First we point out the proof is not exactly straightforward by adding an additional radii for the $b_1(\cdot)$ coefficient since the third term in the above expression of $\sigma^2_i$ also involves $\sigma_{i-1}^2$ which has evolved differently for $\kappa=(\mu(\cdot), a_1(\cdot), b_1(\cdot))$ and $\kappa=(\mu_0(\cdot), a_{10}(\cdot), b_{10}(\cdot))$. 
We begin by estimating the difference of the third term

\begin{eqnarray}\label{eq:blambda1}
|b_{1}(i/n)\sigma^2_{i-1}-b_{01}(i/n)\sigma^2_{0,i-1}|\leq \sigma_{0,i-1}^2\|b_{1}-b_{01}\|_{\infty}+b_{1}(i/n)|\sigma^2_{i-1}-\sigma^2_{0,i-1}|.
\end{eqnarray}

\noindent This leads to the recursion
\begin{align*}
    |\sigma_{i}^2-\sigma_{0i}^2|&\leq\|\mu-\mu_{0}\|_{\infty} + X_{i-1}^2\|a_{1}-a_{01}\|_{\infty}+(1-\frac{M_\mu}{M_X})|\sigma_{i-1}^2-\sigma_{0,i-1}^2| + \sigma_{0,i-1}^2\|b_{1}-b_{01}\|_{\infty}
\end{align*}

\ignore{**Sayar
Continuing this recursion we have 
\begin{align*}
    |\sigma_{i}^2-\sigma_{0i}^2|&\leq  \frac{M_{X}}{M_{\mu}}\|\mu-\mu_{0}\|_{\infty} + \sum_{s=0}^{t-1}(1- \frac{M_{\mu}}{M_X})^{s}\{X_{i-1-s}^2\|a_{1}-a_{01}\|_{\infty} +\sigma_{0,n-1-s}^2\|b_{1}-b_{01}\|_{\infty}\}
\end{align*}}

We have
\begin{align*}
    &\sum_{i=1}^{n-1}\frac{M_{\mu}}{M_{X}}|\sigma_{i}^2-\sigma_{0i}^2|+|\sigma_{n}^2-\sigma_{0i}^2|\\&\quad\leq n\|\mu-\mu_{0}\|_{\infty} + \sum_{i}X_{i-1}^2\|a_{1}-a_{01}\|_{\infty} +(1-\frac{M_{\mu}}{M_{X}})|\sigma_{0}^2-\sigma_{00}^2|+ \sum_i\sigma_{0,i-1}^2|b_{11}-b_{01}|_{\infty}
\end{align*}

As $M_{\mu}<M_X$,
\begin{align*}
    \sum_{i=1}^n|\sigma_{i}^2-\sigma_{0i}^2|\leq & \frac{M_X}{M_{\mu}}\LARGE\{n\|\mu_1-\mu_{0}\|_{\infty} + \sum_{i}X_{i-1}^2\|a_{1}-a_{01}\|_{\infty} \\&\quad+(1-\frac{M_{\mu}}{M_{X}})|\sigma_{0}^2-\sigma_{00}^2|+ \sum_i\sigma_{0,i-1}^2\|b_{1}-b_{01}\|_{\infty}\LARGE\}.
\end{align*}
which will be used for bounding the terms $I$ and $II$ in (\ref{eq:IandII}). Towards bounding the first term $I$, we have, along the lines of (\ref{eq:Rfirstterm}),
\begin{align*}
    \sum_{i=1}^n\IE|\frac{\sigma_{i}^2-\sigma_{0i}^2}{\sigma_{*i}^2}|\leq & \frac{M_X}{\rho M_{\mu}}\LARGE\{n\|\mu_1-\mu_{0}\|_{\infty} + nM_X\|a_{1}-a_{01}\|_{\infty}
    \\&\quad+(1-\frac{M_{\mu}}{M_{X}})|\sigma_{0}^2-\sigma_{00}^2|+ nM_X\|b_{1}-b_{01}\|_{\infty}\LARGE\}.
\end{align*}
For bounding the second term $II$, however, we first look at the summand as follows
\begin{eqnarray}\label{eq:garchcalc}
\IE(X_i^2(1/\sigma_{0i}^2- 1/\sigma_i^2)\})\leq \IE(\IE(\frac{X_i^2}{\sigma_{0i}^2}\frac{|\sigma_{i}^2-\sigma_{0i}^2|}{\sigma_i^2}|X_{i-1})) &\leq& \IE(\frac{|\sigma_{i}^2-\sigma_{0i}^2|}{\sigma_i^2}\IE(\frac{X_i^2}{\sigma_{0i}^2}|X_{i-1})) \\
&\leq& \frac{1}{\rho}\IE(|\sigma_{i}^2-\sigma_{0i}^2|) \nonumber
\end{eqnarray}

\noindent Taking sum over $t$, we have
\begin{align*}
    &\frac{1}{\rho}\sum_i\IE(|\sigma_{i}^2-\sigma_{0i}^2|)=\frac{1}{\rho}\IE(\sum_i|\sigma_{i}^2-\sigma_{0i}^2|)\\
    \quad&\leq\frac{M_X}{\rho M_{\mu}}\LARGE\{n\|\mu_1-\mu_{0}\|_{\infty} + nM_X\|a_{1}-a_{01}\|_{\infty} +n(1-\frac{M_{\mu}}{M_{X}})|\sigma_{0}^2-\sigma_{00}^2|+ nM_X\|b_{1}-b_{01}\|_{\infty}\LARGE\}
\end{align*}

\noindent Combining the bounds for $I$ and $II$ we arrive at 
\begin{align}
   \frac{1}{n}\IE(R)\lesssim \|\mu-\mu_0\|_{\infty}+\|a_1-a_{01}\|_{\infty}+\|a_1-a_{01}\|_{\infty}+\|\sigma_0-\sigma_{00}\|_{\infty}.\label{eq:KLineq2}
\end{align}
In the next part we modify the construction of the sieve to suit the tvGARCH structure.

\subsubsection{Construction of exponentially consistent tests}
\noindent For the part where we construct exponentially consistent test, note that the only challenge that remains for GARCH processes is to obtain a claim as Claim \ref{claim1}. Towards that, we use a new definition of sieve in the light of (\ref{eq:sieve}) as following:
\begin{align}\label{eq:sievegarch}
    W_n&=\{K_1,K_2,K_3,\alpha,\gamma_1, \gamma_2:K_1\leq K_{1n},K_2\leq K_{2n}, K_3\leq K_{3n},\|\alpha\|_{\infty}\leq A_{n}, \min(\alpha,\gamma_1,\gamma_2)>\rho_n,\nonumber\\&\qquad \max{\gamma_1}+\max{\gamma_2}\leq 1-A_n/B_n,\sigma_0^2\leq B_n\},
\end{align}
where $A_n, B_n$ are at least polynomial in $n$ and $\sigma_0^2$ is the standard deviation of $X_0^2$ and $K_n=\max\{K_{1n},K_{2n}\}$. We take $\rho_n \asymp n^{-a}$ with $a<1$, $A_n=B_n(1-n^{1/n}\rho_n), B_n \asymp n^{a_2}$ for sufficiently large $n$ such that $n^{1/n}\rho_n<1$.


\ignore{
************************************************************************************

Note that, in the light of assumption B.2, thanks to the lower bound we can have a new claim as follows 
}
\begin{claim}\label{claim2}
 $\IE_{\kappa_1}(f/f_1)^2$ is bounded for every $\kappa$ such that
\begin{eqnarray}\label{eq:condition2}
r_i\leq \frac{\rho_n}{5n^2}\text{  for }1\leq i \leq 4
\end{eqnarray}
where $\|\mu-\mu_1\|_{\infty}=r_1, \|a-a_1\|_{\infty}=r_2,\|b-b_1\|_{\infty}=r_3, \|\sigma_0^2-\sigma_{01}^2\|=r_4$.
\end{claim}
\begin{proof}

Within the sieve again we use a variant of above inequality. 

\ignore{
**** (Not needed?) 
Note that within the sieve $\IE(X_i^2)\leq B_n$. 

\begin{align*}
    \sum_{i=1}^n|\sigma_{1t}^2-\sigma_{i}^2|\leq & \frac{B_n}{A_n}\LARGE\{n\|\mu_1-\mu\|_{\infty} + \sum_{i}X_{i-1}^2\|a_{11}-a_{1}\|_{\infty} \\&\quad+(1-\frac{A_n}{B_n})|\sigma_{10}^2-\sigma_{0}^2|+ n\max_i\sigma_{i-1}^2|b_{11}-b_{1}|_{\infty}\LARGE\}.
\end{align*}
****
}

\begin{align*}
    \frac{|\sigma_{i}^2-\sigma_{1i}^2|}{\sigma_{i}^2}&\leq\frac{1}{\rho_n}\|\mu-\mu_{1}\|_{\infty} + \frac{1}{\rho_n}\|a_{1}-a_{11}\|_{\infty}+\frac{1-A_n/B_n}{\rho_n}\frac{|\sigma_{i-1}^2-\sigma_{1,n-1}^2|}{\sigma_{i-1}^2} + \frac{1}{\rho_n}\|b_{1}-b_{11}\|_{\infty}
\end{align*}

\noindent By recursion, 
\begin{align}
    \frac{|\sigma_{i}^2-\sigma_{1i}^2|}{\sigma_{i}^2}\leq\frac{Q_n^i-1}{(Q_n-1)\rho_n}\left[\|\mu-\mu_{1}\|_{\infty}+\|a_{1}-a_{11}\|_{\infty}+\|b_{1}-b_{11}\|_{\infty}\right] + \frac{Q_n^{t-1}}{\rho_n}|\sigma_{0}^2-\sigma_{01}^2|,
\end{align}
where $Q_n=\frac{1-A_n/B_n}{\rho_n}>1$. The RHS is increasing in $i$ and thus we only need to find a bound for $i=n$.  $A_n,B_n$ and $\rho_n$ are chosen in such a way that $Q_n  \asymp n^{1/n}$. Based on that $r_1$, $r_2, r_3$ and $r_4$ can be chosen. We also choose $A_n=B_n(1-n^{1/n}\rho_n)$. Note that, with how we choose $\rho_n=n^{-a}$ above, $n^{1/n}\rho_n<1$ for $n>1/a$ which means for large $n$, such choices of $A_n, B_n$ are valid. Finally we make the choices for radii as $r_i\lesssim \rho_n/n^2$ for $1\leq i \leq 4$.

We now conclude the proof of the claim.  Due to the radii choices above and the definition of sieve at (\ref{eq:sievegarch}), we can ensure 
$$\frac{\sigma^2-\sigma_1^2}{\sigma_1^2} \leq  \frac{4}{5n}$$ 
where notation similar to Claim \ref{claim1} is used. This implies $(n+4/5)\sigma_1^2-n\sigma^2>0$, then
\begin{align}
    \int_{x}\frac{\{f(X_i=x|X_{i-1})\}^{n}}{\{f_1(X_i=x|X_{i-1})\}^n}f(X_i=x|X_{i-1})dx &\leq \sqrt{5} \left( \frac{\sigma_1^2}{\sigma^2} \right)^{n/2}\label{eq:tobdexp2}.
\end{align}
along similar lines of (\ref{eq:tobdexp}). Note that, in the light of (\ref{eq:condition2}), for large $n$
$$\left(\frac{\sigma_1^2}{\sigma^2} \right)^{n/2}=\left(1+\frac{\sigma_1^2-\sigma^2}{\sigma^2}\right)^{n/2} \leq \left(1+\frac{4}{5n}\right)^{\frac{n}{2}}\approx e^{2/5}.$$
\end{proof}

\noindent The equivalence of Reyni divergence and Hellinger is exactly same as shown in Subsubsection \ref{ssc:equiv}.

\subsection{Proof of Theorem 3}
We prove Theorem \ref{thm:igarch} along very similar lines as outlined for Theorem \ref{thm:garch}. The only difference pops up in the KL difference step where existence of second moment is necessary but unfortunately for the tviGARCH case, the variance is infinite. Thus we handle the KL bound in a different way. 

\noindent  Our goal is to obtain a deterministic bound on $\sum_{i=1}^n \frac{|\sigma_i^2-\sigma_{i0}^2|}{\sigma_i^2}$. Denoting $a_1(i/n)$ by $a_{1i}$ we have
\begin{eqnarray}
|\sigma_i^2-\sigma_{i0}^2|&\leq & \sigma_{i-1}^2(a_{1i0}-a_{1i}) +(1-a_{1i0}) |\sigma_{i-1}^2-\sigma_{i-1,0}^2|+\|\mu-\mu_{0}\|_{\infty}+(a_{1i}-a_{1i0})X_{i-1}^2 \nonumber \\
& \leq& \frac{1}{\rho} \|\mu-\mu_{0}\|_{\infty}+\|a_1-a_{10}\|(\sigma_{i-1}^2+(1-a_{1i0})\sigma_{i-2}^2+(1-a_{1i0})(1-a_{1,i-1,0})\sigma_{i-3}^2+\cdots)\nonumber \\
&& +\|a_{1}-a_{10}\|_{\infty}(X_{i-1}^2+(1-a_{1i0})X_{i-2}^2+(1-a_{1i0})(1-a_{1,i-1,0})X_{i-3}^2+\cdots)\nonumber \\
\frac{|\sigma_i^2-\sigma_{i0}^2|}{\sigma_i^2}& \leq& \frac{1}{\rho\sigma_i^2} \|\mu-\mu_{0}\|_{\infty}+\frac{1}{\sigma_i^2}(1-a_{1i0})\cdots(1-a_{110})|\sigma_0^2-\sigma_{00}^2|\nonumber \\
&& \quad +\frac{\|a_1-a_{10}\|}{1-\rho} (1+x_i+..+x_i \cdots x_2)+\frac{\|a_1-a_{10}\|}{\rho}(1+x_i+..+x_i \cdots x_2)\nonumber \\
&\leq& \frac{1}{\rho^2} \|\mu-\mu_{0}\|_{\infty}+\frac{1}{\rho}(1-a_{1i0}\cdots(1-a_{110})|\sigma_0^2-\sigma_{00}^2|\nonumber \\ && \quad+\frac{\|a_1-a_{10}\|}{\rho(1-\rho)}(1+x_i+\cdots+x_i \cdots x_2)\nonumber \\ 
&\leq& \frac {\|\mu-\mu_{0}\|_{\infty}}{\rho^2}+\frac{(1-\rho)^i|\sigma_0^2-\sigma_{00}^2|}{\rho}+\frac{\|a_1-a_{10}\|}{\rho(1-\rho)}(1+x_i+..+x_i \cdots x_2) \label{eq:igarchcalc}
\end{eqnarray}
where $x_i= (1-a_{1i0})/(1-a_{1i})$. In the first inequality of the above derivation, we have used $\sigma_i^2$ satisfies the following 
\begin{eqnarray}\label{eq:sigmabound}
\sigma_i^2 \geq \rho,\sigma_i^2&\geq& (1-a_{1i})\sigma_{i-1}^2,\sigma_i^2\geq (1-a_{1i})(1-a_{1,i-1})\sigma_{i-2}^2 ,\cdots, \nonumber \\
\sigma_i^2 &\geq& a_{1i}X_{i-1}^2, \sigma_i^2 \geq (1-a_{1i})a_{1,i-1}X_{i-2}^2, \cdots
\end{eqnarray}
and $\inf_x(\mu(x))>\rho_{\mu},\sup_xa_1(x)<1-\rho$ and $\inf_x a_1(x)>\rho$ thanks to the closeness of $(\mu,a)$ with $(\mu_0,a_{10})$. For the third term, we have bounds on $x_i$ as follows
$x_i\leq 1+ (\|a_1-a_{10}\|_{\infty})/\rho$. Consequently,  
\begin{eqnarray*}
\sum_{i=1}^n \frac{\|a_1-a_{10}\|}{\rho(1-\rho)}(1+x_i+\cdots+x_i \cdots x_2) &\leq& \frac{\|a_1-a_{10}\|}{\rho(1-\rho)}\sum_{i=1}^n\sum_{j=0}^{i-1}\left( 1+ \frac{\|a_1-a_{10}\|_{\infty}}{\rho}\right)^j \\ &=& \frac{1}{(1-\rho)}\sum_{i=1}^n\left( 1+ \frac{\|a_1-a_{10}\|_{\infty}}{\rho}\right)^i- \frac{n}{1-\rho} \\
&=& \frac{\|a_1-a_{10}\|_{\infty}+\rho}{\rho(1-\rho)} \times\frac{\left( 1+ \frac{\|a_1-a_{10}\|_{\infty}}{\rho}\right)^{n}-1}{\frac{\|a_1-a_{10}\|_{\infty}}{\rho}}-\frac{n}{1-\rho}\\
&\approx& \frac{\|a_1-a_{10}\|_{\infty}+\rho}{\rho(1-\rho)} \times \frac{\left( 1+\frac{n\|a_1-a_{10}\|_{\infty}}{\rho}-1\right)}{\frac{\|a_1-a_{10}\|_{\infty}}{\rho}} -\frac{n}{1-\rho} \\
&=& \frac{n\|a_1-a_{10}\|_{\infty}}{\rho(1-\rho)}
\end{eqnarray*}
by a Taylor series approximation of the binomial term in the third line since $\|a_1-a_{10}\|_{\infty}$ is small. This along with (\ref{eq:igarchcalc}) leads to the following bound 
$$\sum_{i=1}^n\frac{|\sigma_i^2-\sigma_{i0}^2|}{\sigma_i^2}\leq \frac{n}{\rho^2}\|\mu-\mu_{0}\|_{\infty}+\frac{1}{\rho^2}|\sigma_0^2-\sigma_{00}^2|+\frac{n\|a_1-a_{10}\|_{\infty}}{\rho(1-\rho)}
\lesssim n\|\mu-\mu_0\|+|\sigma_0^2-\sigma_{00}^2|+n\| a_1-a_{10}\|_{\infty}.$$

\noindent Note that for bounding $I$ in the KL term from (\ref{eq:IandII}), $\sigma_{i}^{*2}$ satisfy the exact same bounds as (\ref{eq:sigmabound}). For term $II$ as well the bounds follow along the lines of (\ref{eq:garchcalc}).

We consider following sieve for iGARCH:
\begin{align}\label{eq:sieveigarch}
    W_n&=\{K_1,K_2,\alpha,\gamma_1,:K_1\leq K_{1n},K_2\leq K_{2n},\|\alpha\|_{\infty}\leq A_{n}, \min(\alpha,\gamma_1)>\rho_n,\nonumber\\&\qquad \max{\gamma_1}\leq 1-\rho_n,\sigma_0^2\leq B_n\},
\end{align}
Within the sieve,
$$\sum_{i=1}^n\frac{|\sigma_i^2-\sigma_{i1}^2|}{\sigma_i^2}\leq \frac{n}{\rho_n^2}\|\mu-\mu_{1}\|_{\infty}+\frac{1}{\rho_n^2}|\sigma_0^2-\sigma_{01}^2|+\frac{n\|a_1-a_{11}\|_{\infty}}{\rho_n(1-\rho_n)}.$$
Along the line of Claim \ref{claim1} and Claim \ref{claim2}, the rest of the proof goes through if the radii satisfy $\|\mu-\mu_{1}\|_{\infty}=r_1\leq \frac{\rho_n^2}{8n}$, $|\sigma_0^2-\sigma_{01}^2|=r_2 \leq \rho_n^2/8$ and $\|a_1-a_{11}\|_{\infty}=r_3 \leq \frac{\rho_n(1-\rho_n)}{8n}$.

\ignore{\subsection{Proof of Theorem 3}

\begin{align*}
    \sigma_{i}^2-\sigma_{0i}^2 = \mu(i/n)-\mu_{0}(i/n) + (X_{i-1}^2-\sigma_{0,i-1}^2)(a_{1}(i/n)-a_{01}(i/n))+(1-a_{1}(i/n))(\sigma_{i-1}^2-\sigma_{0,i-1}^2)
\end{align*}
which means
\begin{align*}
    \IE_{\psi_0}(\sigma_{i}^2-\sigma_{0i}^2) = \mu(i/n)-\mu_{0}(i/n)-(1-a_{01}(i/n))(\sigma_{i-1}^2-\sigma_{0,i-1}^2) + (a_{01}(i/n)-a_{1}(i/n))(\sigma_{0,i-1}^2-\sigma_{i-1}^2)
\end{align*}
Since $(1-\rho)>a_{01}(i/n)>\rho$ for all $i$, we have
\begin{align*}
    \IE_{\psi_0}(|\sigma_{i}^2-\sigma_{0i}^2|) \leq \|\mu-\mu_{0}\|+\rho\IE_{\psi_0}|\sigma_{0,i-1}^2-\sigma_{i-1}^2| + \|a_{01}-a_1\|_{\infty}\IE_{\psi_0}|\sigma_{0,i-1}^2-\sigma_{i-1}^2|
\end{align*}
 Let $\rho+\|a_{01}-a_1\|_{\infty}=C_1<1$ for $\|a_{01}-a_1\|_{\infty}<\rho$. Then by recurssion,
 \begin{align*}
    \IE_{\psi_0}(|\sigma_{i}^2-\sigma_{0i}^2|) \leq \|\mu-\mu_{0}\|\frac{1}{1-C_1}+C_1^t|\sigma_{0,0}^2-\sigma_{0}^2|
\end{align*}

 \subsubsection{Test construction}
\begin{align*}
    \frac{|\sigma_{i}^2-\sigma_{1t}^2|}{\sigma_{i}^2}=\frac{1}{\sigma_{i}^2}(\mu(i/n)-\mu_{1}(i/n) + \frac{1}{\rho_n}\|a_{1}-a_{11}\|_{\infty}+\frac{1-A_n/B_n}{\rho_n}\frac{|\sigma_{i-1}^2-\sigma_{1,n-1}^2|}{\sigma_{i-1}^2} + \frac{1}{\rho_n}\|b_{1}-b_{11}\|_{\infty}
\end{align*} 
}

\newpage

\begin{figure}
    \centering
    \includegraphics[width=100mm]{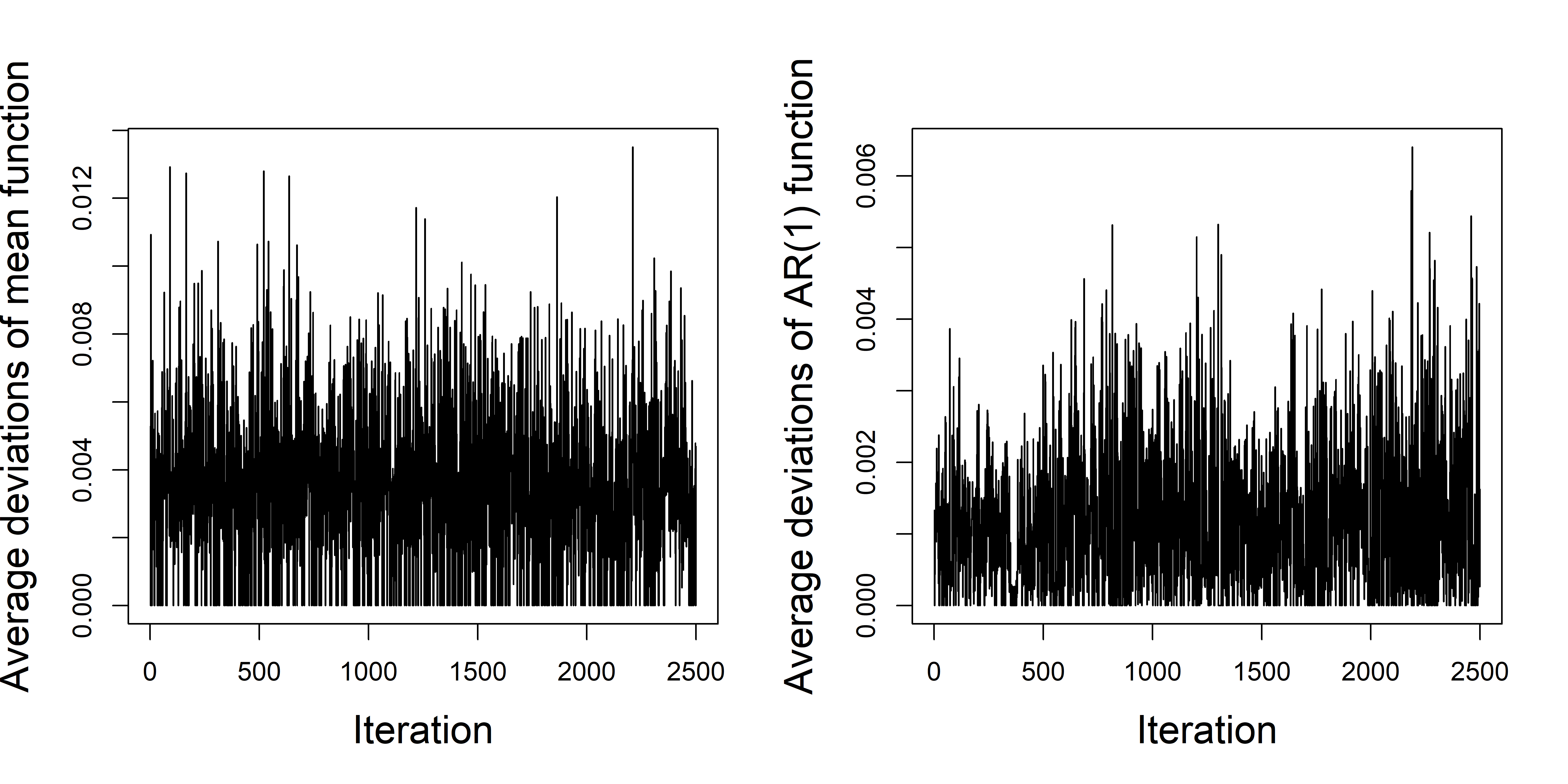}
    \caption{Trace plots of average $L_{2}$ deviations across the MCMC chain, thinned by 2 for the mean function $\frac{1}{1000}\left[\sum_{j=1}^{1000}\overline{(\mu(j/1000)^t-\mu(j/1000)^{t+1})^2}\right]^{1/2}$ and for the AR(1) function $\frac{1}{1000}\left[\sum_{j=1}^{1000}\overline{(a_1(j/n)^t-a_1(j/n)^{t+1})^2}\right]^{1/2}$. Here $\mu(j/1000)^t$ and $a_1(j/n)^t$ are the $t$-$th$ post burn samples of $\mu(j/1000)$ and $a_1(j/n)$ respectively for ARCH model.}
    \label{fig:trARCH}
\end{figure}

\begin{figure}
    \centering
    \includegraphics[width=100mm]{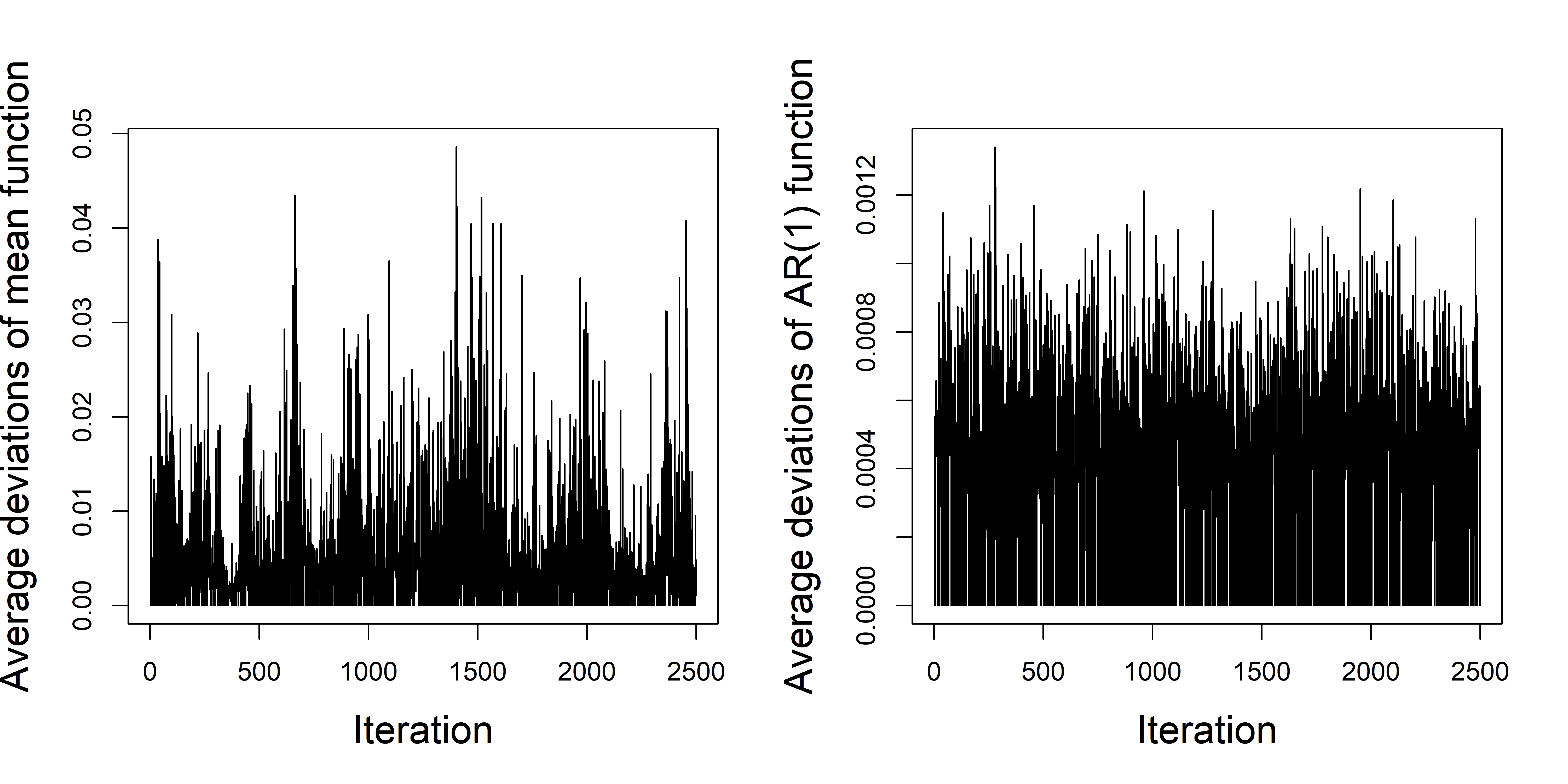}
    \caption{Trace plots of average $L_{2}$ deviations across the MCMC chain, thinned by 2 for the mean function $\frac{1}{1000}\left[\sum_{j=1}^{1000}\overline{(\mu(j/1000)^t-\mu(j/1000)^{t+1})^2}\right]^{1/2}$ and for the AR(1) function $\frac{1}{1000}\left[\sum_{j=1}^{1000}\overline{(a_1(j/n)^t-a_1(j/n)^{t+1})^2}\right]^{1/2}$. Here $\mu(j/1000)^t$ and $a_1(j/n)^t$ are the $t$-$th$ post burn samples of $\mu(j/1000)$ and $a_1(j/n)$ respectively for IGARCH model.}
    \label{fig:trIGARCH}
\end{figure}

\begin{figure}
    \centering
    \includegraphics[width=120mm]{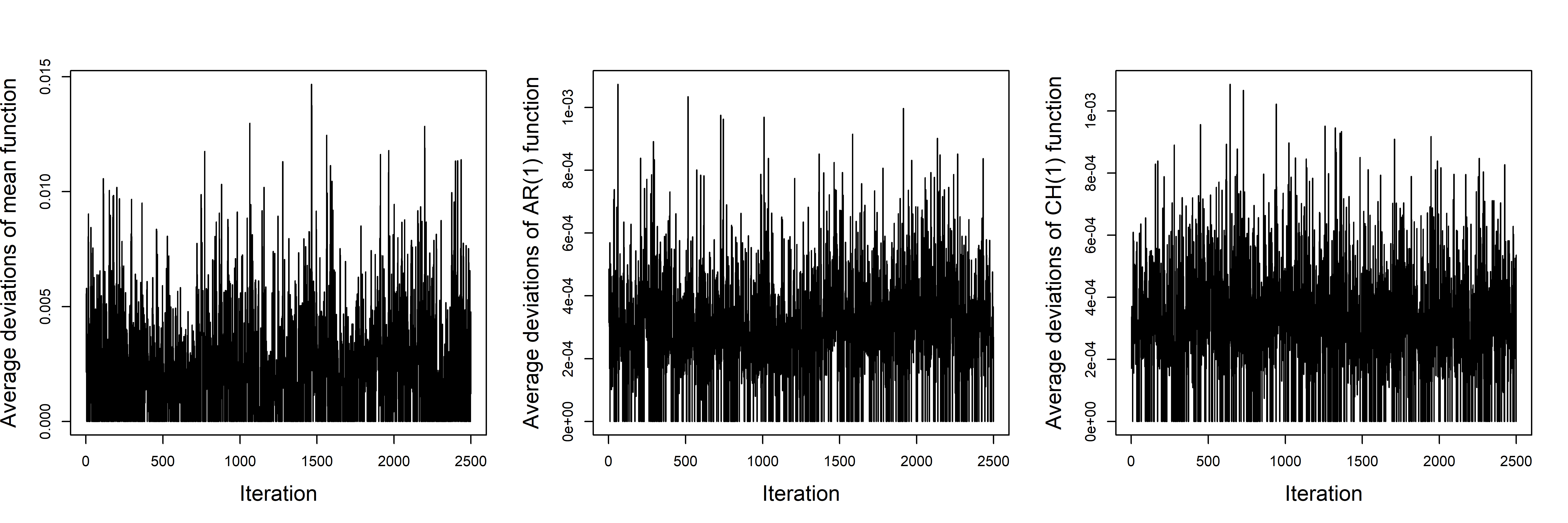}
    \caption{Trace plots of average $L_{2}$ deviations across the MCMC chain, thinned by 2 for the mean function $\frac{1}{1000}\left[\sum_{j=1}^{1000}\overline{(\mu(j/1000)^t-\mu(j/1000)^{t+1})^2}\right]^{1/2}$, for the AR(1) function $\frac{1}{1000}\left[\sum_{j=1}^{1000}\overline{(a_1(j/n)^t-a_1(j/n)^{t+1})^2}\right]^{1/2}$, and for the CH(1) function $\frac{1}{1000}\left[\sum_{j=1}^{1000}\overline{(b_1(j/n)^t-b_1(j/n)^{t+1})^2}\right]^{1/2}$. Here $\mu(j/1000)^t$, $a_1(j/n)^t$, and $b_1(j/n)^t$ are the $t$-$th$ post burn samples of $\mu(j/1000)$ , $a_1(j/n)$, and $b_1(j/n)$ respectively for GARCH model.}
    \label{fig:trGARCH}
\end{figure}

\end{document}